\documentclass[a4paper]{amsart}
\usepackage{amssymb}
\usepackage{amscd}
\usepackage{enumitem}
\setlist{font=\normalfont}
\usepackage{graphicx}
\usepackage{color}
\usepackage{hyperref}
\usepackage{nccmath}
\usepackage{hyphenat}
\usepackage{makecell}
\usepackage{float}

\theoremstyle{plain}
\newtheorem{theorem}{Theorem}
\newtheorem{lemma}[theorem]{Lemma}
\newtheorem{corollary}[theorem]{Corollary}
\newtheorem{proposition}[theorem]{Proposition}
\newtheorem{fact}[theorem]{Fact}

\theoremstyle{definition}
\newtheorem{definition}[theorem]{Definition}
\newtheorem{example}[theorem]{Example}

\newtheorem{remark}[theorem]{Remark}

\numberwithin{equation}{section}
\numberwithin{theorem}{section}

\renewcommand{\a}{\mathbf{a}}
\newcommand{\x}{\mathbf{x}}
\newcommand{\y}{\mathbf{y}}
\newcommand{\s}{\mathbf{s}}
\newcommand{\teta}{\vartheta}
\newcommand{\2}{\mathbf{2}}
\newcommand{\N}{\mathbb{N}}

\renewcommand{\O}{\mathcal{O}}
\newcommand{\OA}{\mathcal{O}_A}
\newcommand{\OAn}{\mathcal{O}_A^{(n)}}
\newcommand{\OL}{\mathcal{O}_L}

\newcommand{\OS}{\mathcal{O}_S}

\newcommand{\down}[1]{{\downarrow}\,#1}
\newcommand{\up}[1]{{\uparrow}\,#1}

\newcommand{\Hom}{\operatorname{Hom}}
\newcommand{\Homj}{\Hom_{\lor}}
\newcommand{\Homm}{\Hom_{\land}}
\newcommand{\Con}{\operatorname{Con}}

\begin{document}

\title[On Centralizers of finite lattices and semilattices]{On Centralizers of finite lattices and semilattices}

\author[E. T\'{o}th]{Endre T\'{o}th}
\address[E. T\'{o}th]{Bolyai Institute, University of Szeged, Aradi v\'{e}rtan\'{u}k tere 1, H--6720 Szeged, Hungary}
\email{tothendre@math.u-szeged.hu}
\author[T. Waldhauser]{Tam\'as Waldhauser}
\address[T. Waldhauser]{Bolyai Institute, University of Szeged, Aradi v\'{e}rtan\'{u}k tere 1, H--6720 Szeged, Hungary}
\email{twaldha@math.u-szeged.hu}

\begin{abstract}
We study centralizer clones of finite lattices and semilattices.
For semilattices, we give two characterizations of the centralizer and also derive formulas for the number of operations of a given essential arity in the centralizer.
We also characterize operations in the centralizer clone of a distributive lattice, and we prove that the essential arity of operations in the centralizer is bounded for every finite (possibly nondistributive) lattice.
Using these results, we present a simple derivation for the centralizers of clones of Boolean functions.

\end{abstract}



\maketitle

\section{Introduction}
\begin{figure}[t]
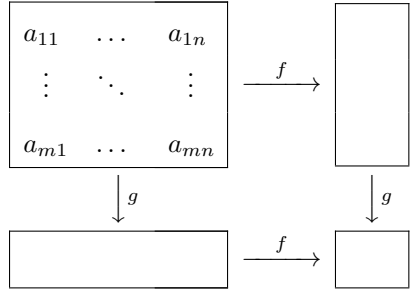

\[
\begin{CD} \renewcommand{\arraystretch}{1.6}\begin{tabular}{|m{0.6cm}m{0.6cm}m{0.6cm}|} \cline{1-3}\cline{3-3} $a_{11}$ & $\dots$ & $a_{1n}$ \\ \hfil$\vdots$ & $\ddots$ & \hfil$\vdots$\\ $a_{m1}$ & $\dots$ & $a_{mn}$ \\ \cline{1-3}\cline{3-3} \end{tabular} @>f>> \renewcommand{\arraystretch}{1.66}\begin{tabular}{m{0.6cm}} \cline{1-1} \multicolumn{1}{|m{0.6cm}|}{} \\ \multicolumn{1}{|c|}{} \\ \multicolumn{1}{|c|}{} \\ \cline{1-1} \end{tabular} \\ @VVgV @VVgV\\ \renewcommand{\arraystretch}{1.66}\begin{tabular}{|m{0.6cm}m{0.6cm}m{0.6cm}|} \cline{1-3}\cline{3-3} & & \\ \cline{1-3}\cline{3-3} \end{tabular} @>f>> \renewcommand{\arraystretch}{1.66}\begin{tabular}{|m{0.6cm}|} \cline{1-1} \\ \cline{1-1} \end{tabular} \end{CD}
\]
\caption{Commutation of $f$ and $g$.} \label{fig:comm}
\end{figure}
We say that the operations $f\colon A^n \to A$ and $g\colon A^m \to A$ \emph{commute } (notation: $f\perp g$) if%
\begin{multline*}
g\bigl(f(a_{11},a_{12},\dotsc,a_{1n}),f(a_{21},a_{22},\dotsc,a_{2n}),\dotsc,f(a_{m1},a_{m2},\dotsc
,a_{mn})\bigr)\\
=f\bigl(g(a_{11},a_{21},\dotsc,a_{m1}),g(a_{12},a_{22},\dotsc,a_{m2}),\dotsc,g(a_{1n},a_{2n},\dotsc
,a_{mn})\bigr)\\
\end{multline*}
holds for all $a_{ij}\in A~(1\leq i\leq m,1\leq j\leq n)$. This can be visualized as follows: for every $m\times n$ matrix $Q=(a_{ij})$, first applying $f$ to the rows of $Q$ and then applying $g$ to the resulting column vector yields the same result as first applying $g$ to the columns of $Q$ and then applying $f$ to the resulting row vector (see Figure~\ref{fig:comm}).
In other words, $f \perp g$ holds if and only if $g$ is a homomorphism from $(A;f)^m$ to $(A;f)$.
In particular, if $g$ is a unary operation (i.e., $m=1$), then $g$ commutes with $f$ if and only if $g$ is an endomorphism of the algebra $(A;f)$.

Commutation of operations is one of the many areas where Ivo Rosenberg had outstanding achievements, and we would like to contribute to this memorial issue by some results that were inspired by one of his last papers \cite{MachidaRosenberg} on this topic. 
We study operations commuting with (semi)lattice operations, and we aim for concrete descriptions of these operations that allow us to classify (and in some cases also count) them according to the number of their essential variables.
As a ``byproduct", we also obtain a simple proof for Kuznetsov's description \cite{Kuz} of primitive positive clones on the two-element set.

Let us introduce some notions and notation that will be used throughout the paper. 
An $n$-ary \emph{operation} on a set $A$ (which will always assumed to be finite) is a map $f\colon A^n \to A$.
The set of all $n$-ary operations on $A$ is denoted by $\OAn$, and the set of all finitary operations on $A$ is $\OA=\bigcup_{n\geq0}\OAn$.
We say that the $i$-th variable of $f\in\OAn$ is \emph{essential} (or that $f$ \emph{depends} on its $i$-th variable) if there exist tuples $\a,\a' \in A^n$ differing only in their $i$-th component such that $f(\a)\neq f(\a')$.
The number of essential variables of $f$ is called the \emph{essential arity} of $f$.
To simplify notation, we often assume that operations do not have inessential variables, thus we say that $f$ is an essentially $n$-ary operation on $A$ if $f\in\OAn$ and $f$ depends on all of its variables.
 
A set $C$ of operations on $A$ is a \emph{clone} (notation: $C\leq\OA$) if $C$ is closed under composition and contains the \emph{projections} $(x_1,\dots,x_n)\mapsto x_i$ for all $1\leq i\leq n$. 
The least clone containing a set $F\subseteq\OA$  is called the clone \emph{generated} by $F$, and it is denoted by $[F]$. 
The \emph{centralizer} of $F$ is the set $F^\ast$ of all operations commuting with each member of $F$:
\[
F^\ast = \{ g\in\OA : f\perp g \text{ for all } f \in F\}.
\]
\noindent It is easy to verify that $F^\ast$ is a clone (even if $F$ is not), and $F^\ast=[F]^\ast$. This means that the centralizer of the clone of term operations of an algebra $(A;F)$  is $F^\ast$, which we call simply the centralizer of the algebra $(A;F)$.

Although there are countably infinitely many clones on the two-element set \cite{Post} and uncountably many clones on sets with at least three elements \cite{JM}, only finitely many clones are of the form $F^\ast$ on a finite set \cite{BW}. 
These are the so-called \emph{primitive positive clones}; the complete list of primitive positive clones is known only for $|A|\leq 3$ \cite{Kuz, Dan}.

In Section~\ref{sect:semilattices} we present two different characterizations of the essentially $n$-ary members of the centralizer of a finite semilattice (one of them is a slight variation of a result of Larose \cite{Larose}). 
We use these to  give a general formula for the number of essentially $n$-ary operations commuting with the join operation of a finite lattice, and we illustrate this with the example of finite chains. 
This is a generalization of one of the results of \cite{MachidaRosenberg}, where this counting problem was solved for the three-element chain. 
We study operations commuting with both the join and the meet operation of a lattice in sections~\ref{sect:distributive} and \ref{sect:lattices}. 
Since the essential arity of operations in the centralizer is bounded for every finite lattice, here we focus on the existence of essentially $n$-ary operations instead of counting them. 
In Section~\ref{sect:distributive} we give an explicit description of the elements of the centralizer of a finite distributive lattice, and then we provide two characterizations of finite distributive lattices having an essentially $n$-ary operation in their centralizers.
In Section~\ref{sect:lattices} we investigate thoroughly all results of the previous section to see which ones remain valid for nondistributive finite lattices. 
As a tool aiding this investigation, we give an upper bound for the essential arity of operations in the centralizer of a finite algebra generating a congruence distributive variety (following an idea of \cite{CGL, CM}).
Finally, in Section~\ref{sect:boole} we describe the centralizer clones of Boolean functions. 
Primitive positive clones on the two-element set were described already by Kuznetsov \cite{Kuz}, and later also in \cite{Miki}, and probably many of the readers of this paper have also computed these by themselves at some point. 
This is not a difficult task using the Post lattice (see Figure~\ref{fig:post} in the appendix), but it involves some case-by-case analysis.
We offer a ``painless" proof that covers all cases by just three general theorems.
Besides presenting the list of the $25$ primitive positive clones, we also give the centralizer of each Boolean clone in Table~\ref{table centralizers on 01}; we hope that it serves as a useful reference for anyone studying centralizer clones.

Some personal remarks from the second author about Ivo Rosenberg: 
For several years, Rosenberg's Theorem on the five types of minimal clones was my daily bread, as I was a doctoral student working on minimal clones under the supervision of B\'{e}la Cs\'{a}k\'{a}ny and \'{A}gnes Szendrei.
I met Ivo Rosenberg only once, and even then we spoke only a few words, but he still helped me a lot at the start of my career.
He handled my very first paper as an editor of Algebra Universalis, and his encouragement meant a lot for me, when, to my horror, the referee found a serious gap in the proof of the main result.
It took a lot of work to fill the gap, but eventually the paper was fixed and published.
Later, Ivo Rosenberg was one of the reviewers of my PhD thesis.
Of course, he did not fly from Montr\'{e}al to Szeged for the defense, but, besides the official report, he sent a friendly hand-written letter to me (quite a curiosity nowadays!), in which he pointed out a few typos in the thesis (he was kind enough not to put them into the report) and expressed his interest in my work.
Time flies quickly, the next generation is already here: now I am writing this paper with my doctoral student Endre T\'{o}th, and we dedicate this contribution to the memory of Ivo Rosenberg.

\section{Centralizers of finite semilattices}\label{sect:semilattices}
In this section we give two different characterizations of the centralizer clone of a
join\hyp{}semilattice $S=(S;\lor)$. Since we are interested in counting the 
essentially $n$-ary operations in the centralizer, we assume that $S$ is finite,
but some of our results are also valid for infinite complete semilattices.

\subsection{Characterizations}
If $S$ is a finite join\hyp{}semilattice then it has a greatest element (denoted by $1$),
and if $S$ also has a least element (denoted by $0$), then there is a meet operation
on $S$ such that $(S;\lor,\land)$ is a lattice.
In the latter case the centralizer clone $[\lor]^\ast$ is generated by
its unary members (i.e., endomorphisms of $(S;\lor)$) together with the join
operation. 
This was proved by B.~Larose \cite{Larose}; in the following theorem 
we reprove this result, and we extend it by providing a unique expression for any 
$f \in [\lor]^\ast$ as a join of endomorphisms, and we also determine the 
necessary and sufficient condition for $f$ to depend on all of its variables.

\begin{theorem} \label{thm Larose}
Let $S=(S; \lor)$ be a finite semilattice with a least element $0$ and 
greatest element $1$. 
An $n$-ary operation $f \in \OS$  belongs to the centralizer $[\lor]^\ast$ if and only if 
there exist unary operations $u_1, \dots, u_n \in [\lor]^\ast$ such that
\[
f(x_1, \dots, x_n) = u_1(x_1) \lor \dots \lor u_n(x_n)
\text{ and } u_1(0)=u_2(0)= \dots = u_n(0).
\]
The above expression for $f$ is unique, and $f$ depends on all of its variables 
if and only if none of the $u_i$ are constant, i.e., 
$u_i(0) \neq u_i(1)$ for all $i \in \{ 1,\dots,n \}$.
\end{theorem}

\begin{proof}
The ``if" part is clear: the join operation commutes with itself, hence if $f$ can be written as 
a composition of $\lor$ with $u_1, \dots, u_n \in [\lor]^\ast$, 
then $f \in [\lor]^\ast$, as $[\lor]^\ast$ is closed under composition.

For the ``only if" part let us assume that $f$ is an $n$-ary operation in $[\lor]^\ast$, and 
define $u_1, \dots, u_n$ by $u_1(x)=f(x, 0, \dots, 0), \dots, u_n(x)=f(0, \dots, 0, x)$. 
Then obviously $u_1(0)=\dots=u_n(0)=f(0, \dots, 0)$; furthermore, $u_1, \dots, u_n \in [\lor]^\ast$, 
as $f$ and the constant $0$ operation belong to $[\lor]^\ast$. 
Since $f$ commutes with $\lor$, applying the definition of commutation to the $n \times n$ matrix
$$\begin{pmatrix}
x_1 & 0 & 0 & \dots & 0 \\
0 & x_2 & 0 & \dots & 0 \\
\vdots & \vdots & \vdots & \ddots & \vdots \\
0 & 0 & 0 & \dots & x_n \\
\end{pmatrix},$$
\noindent
we can conclude that
\begin{align*} f(x_1, \dots, x_n)
&= f\big(x_1 \lor 0 \lor \dots \lor 0,\ \dots,\ 0 \lor \dots \lor 0 \lor x_n\big) \\
&= f(x_1, 0, \dots, 0) \lor \dots \lor f(0, \dots, 0, x_n) \\
&= u_1(x_1) \lor \dots \lor u_n(x_n).
\end{align*}

To prove uniqueness, assume that $f(x_1, \dots, x_n) = u_1(x_1) \lor \dots \lor u_n(x_n)$
and $u_1(0)=u_2(0)= \dots = u_n(0)$. Then we have 
\[
f(x_1,0, \dots, 0)=u_1(x_1) \lor u_2(0) \lor \dots \lor u_n(0) = u_1(x_1) \lor u_1(0) \lor \dots \lor u_1(0) = u_1(x_1),
\]
as $u_1$ is monotone. Thus $u_1$ is indeed uniquely determined by $f$, and
the above equality also shows that if $u_1(0) \neq u_1(1)$, then $f$ depends 
on its first variable: $f(0,0, \dots, 0)=u_1(0)\neq u_1(1)=f(1,0, \dots, 0)$.
The statements about the uniqueness of $u_i$ and about the essentiality of the 
$i$-th variable for $i=2,\dots,n$ can be proved in an analogous way.
\end{proof}

\begin{figure}[t]
    \centering
    \includegraphics[width=0.4\linewidth]{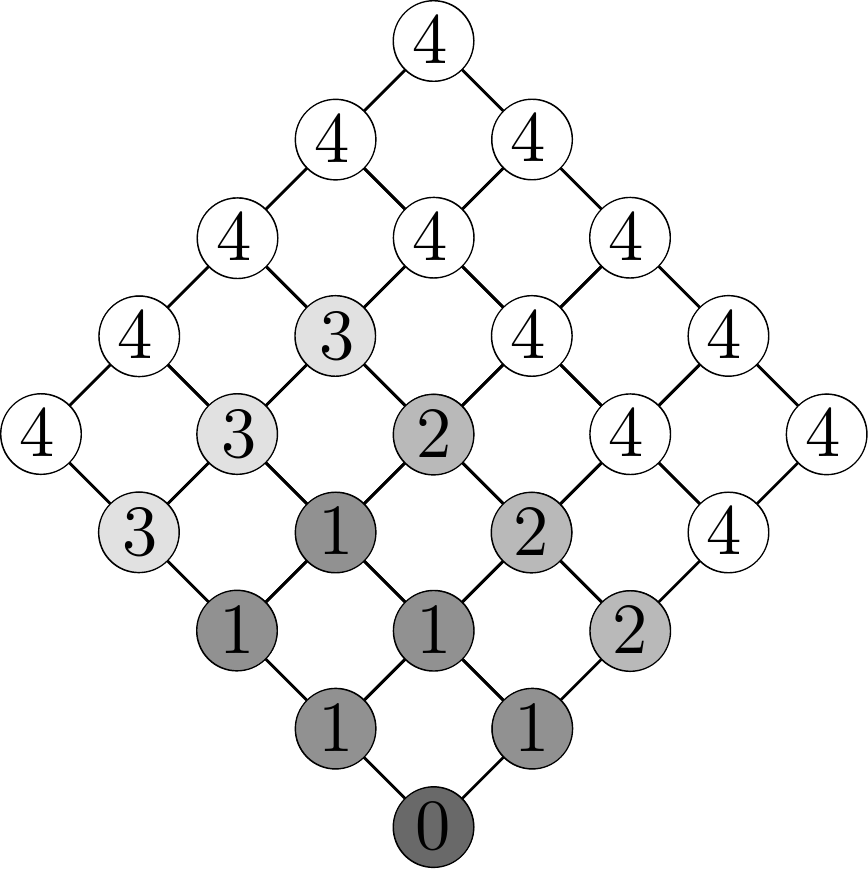}
    \caption{A binary operation in the centralizer of the chain $(\{0,1,2,3,4\};\lor)$}
    \label{fig chain}
\end{figure}

\begin{example}
Let $S=(\{0,1,2,3,4\};\lor)$ be a five-element chain (regarded as a join-semilattice), 
and let the unary operations $u_1$ and $u_2$ be defined by
\begin{align*}
&u_1(0)=0,\ u_1(1)=1,\ u_1(2)=1,\ u_1(3)=3,\ u_1(4)=4;\\
&u_2(0)=0,\ u_2(1)=1,\ u_2(2)=2,\ u_2(3)=4,\ u_2(4)=4.
\end{align*}
Figure~\ref{fig chain} shows the values of the binary operation 
$f(x_1,x_2)=u_1(x_1) \lor u_2(x_2)$. We can see that each of the sets 
$\{ (a_1,a_2) \in S^2 : f(a_1,a_2) \leq b\}\ (b=0,\dots,4)$ is a ``lower rectangle".
We will see later that for every join-semilattice $S$ and every $b \in S$,
the set $\{ \a \in S^n : f(\a) \leq b\}$ has a similar structure for
all $f \in [\lor]^\ast$; moreover, we will characterize the operations
in the centralizer in terms of these ``down-sets".
\end{example}

The following theorem shows that the assumption about $S$ having a least element 
cannot be dropped from Theorem~\ref{thm Larose}: if a finite join\hyp{}semilattice 
does not have a least element, then its centralizer cannot be generated by unary operations and the join operation.

\begin{theorem} \label{thm not generated by unary and join}
The centralizer $[\lor]^\ast$ of a finite semilattice $S=(S;\lor)$ is generated by 
its unary part and the join operation if and only if $S$ has a least element
(i.e., if $S$ is the join\hyp{}reduct of a lattice).
\end{theorem}
\begin{proof}
The ``if" part follows from Theorem~\ref{thm Larose}; for the ``only if" part
assume that $S=(S; \lor)$ is a finite semilattice without a least element.
Then there are distinct minimal elements $a,b \in S$.
We define a binary operation $f$ on $S$ by
\[
f(x_1,x_2)=
\begin{cases} 
   a, & \text{if } (x_1,x_2)=(a,b); \\
   b, & \text{if } (x_1,x_2)=(b,b); \\
   a \lor b, & \text{otherwise}.
\end{cases}
\]
In order to prove that $f$ commutes with the join operation,
we need to verify the following identity:
\begin{equation} \label{eq not generated by unary and join}
f(x_1,x_2) \lor f(y_1,y_2) = f(x_1 \lor y_1, x_2 \lor y_2).
\end{equation}
Since $a$ is a minimal element, the only way of writing $a$ as the join of two elements is $a = a \lor a$. Therefore,
the left hand side of \eqref{eq not generated by unary and join} is $a$
if and only if $f(x_1,x_2) = f(y_1,y_2) = a$, and, by the definition of $f$, 
this holds only for $x_1 = y_1 = a$ and $x_2 = y_2 = b$. The right hand side
of \eqref{eq not generated by unary and join} equals $a$ if and only if
$x_1 \lor y_1 = a$ and $x_2 \lor y_2 = b$, and this is also equivalent to
$x_1 = y_1 = a$ and $x_2 = y_2 = b$. Similarly, both the left hand side and
the right hand side of \eqref{eq not generated by unary and join} take the
value $b$ if and only if $x_1 = y_1 = x_2 = y_2 = b$. For all other
inputs, both sides of \eqref{eq not generated by unary and join} give $a \lor b$.
Thus $f$ belongs to the centralizer $[\lor]^\ast$, indeed.

If a binary operation can be obtained as a composition of the 
join operation and endomorphisms of $S$, then it can be written as 
$u_1(x_1) \lor u_2(x_2)$, where $u_1$ and $u_2$ are endomorphisms. 
(The proof of this fact is a routine term induction; we leave it to the reader.) 
Assume, for contradiction, that our operation $f$ can be expressed
in this form: $f(x_1,x_2)=u_1(x_1) \lor u_2(x_2)$. 
Then we have $a = f(a,b) = u_1(a) \lor u_2(b)$,
and this implies $u_1(a) = u_2(b) = a$. 
On the other hand, $b = f(b,b) = u_1(b) \lor u_2(b) = u_1(b)  \lor a \geq a$, 
which is a contradiction.
\end{proof}

Next we derive another kind of characterization of the centralizer of 
the clone $[\lor]$, which describes, in some sense, the ``distribution" of the values 
of an $n$-ary operation $f \in [\lor]^\ast$ on $S^n$, as illustrated by 
Figure~\ref{fig chain}. For this characterization we will not need the assumption
that $S$ has a least element. Nevertheless, we will consider the
lattice $S_\bot$ obtained from $S$ by adding a new element $\bot$ to the 
bottom of $S$. Thus let $S_\bot = S \cup \{ \bot \}$, where $\bot$ is
an element not contained in $S$, and we define the partial order on
$S_\bot$ so that $\bot < a$ for all $a \in S$, and we keep the original 
ordering on the elements of $S$. Note that we add a new bottom element even 
if $S$ happens to have a least element $0$; in this case $\bot$ is the unique
lower cover of $0$.

If $S$ is a join semilattice, then $f \in [\lor]^\ast$ 
if and only if $f$ is a join\hyp{}homomorphism from $S^n$ to $S$.
This motivates us to consider the set $\Homj(A,B)$ of all join\hyp{}homomorphisms
from $A$ to $B$, where $A$ and $B$ are finite join\hyp{}semilattices.
We use the notation $\Homj^{1}(A,B)$ for the set of all join\hyp{}homomorphisms 
from $A$ to $B$ that preserve the greatest element: 
$\Homj^{1}(A,B):= \{ f \in \Homj(A,B) : f(1_A)=1_B \}$.
Similarly, if $A$ and $B$ have a least element, denoted by 0, then 
let $\Homj^{0}(A,B)$ denote the set 
of join\hyp{}homomorphisms preserving the least element, and let 
$\Homj^{01}(A,B)$ be the set of join\hyp{}homomorphisms preserving both 
boundary elements. Note that the least element of $A_\bot$ and $B_\bot$
is denoted by $\bot$ (not by $0$), therefore in this case we will use
the notation $\Homj^{\bot,1}(A_\bot,B_\bot)$ 
instead of $\Homj^{01}(A_\bot,B_\bot)$.
For meet\hyp{}semilattices $A$ and $B$, the sets $\Homm(A,B)$, $\Homm^{0}(A,B)$, 
etc. are defined analogously.

We need to introduce one more notation: for an element $a$ in a 
partially ordered set ($A;\leq)$, the 
\emph{principal ideal} and the \emph{principal filter} generated by $a$ are 
defined and denoted as follows:
\[
\down{a} = \{ c \in A : c \leq a \}, \quad\quad
\up{a} = \{ c \in A : c \geq a \}.
\]
Observe that in Figure~\ref{fig chain}, the elements labeled by numbers less than
or equal to $b$ form a principal ideal for any $b \in \{0,\dots,4\}$.
This is a special case of the following lemma, which states that 
for all $f \in \Homj(A,B)$ and $b \in B$,
the set $f^{-1}(\down{b})=\{ a \in A : f(a) \leq b\}$ is a principal
ideal whenever it is not empty.

\begin{lemma} \label{lemma inverse image of ideal}
Let $A=(A;\lor)$ and $B=(B;\lor)$ be finite semilattices.
If $f\colon A \to B$ is a homomorphism, then $f^{-1}(\down{b}) \subseteq A$ is 
either empty or a principal ideal for all $b \in B$.
\end{lemma}

\begin{proof}
Assume that $f^{-1}(\down{b})$ is not empty.
Join\hyp{}homomorphisms are monotone, hence it is clear that $f^{-1}(\down{b})$ is 
an ideal, i.e.,
$a_1 \leq a_2 \in f^{-1}(\down{b})$ implies that $a_1 \in f^{-1}(\down{b})$.
Moreover, $f^{-1}(\down{b})$ is closed under joins:
if $a_1,a_2 \in f^{-1}(\down{b})$, then $f(a_1 \lor a_2)=f(a_1) \lor f(a_2) \leq b \lor b =b$,
hence $a_1 \lor a_2 \in f^{-1}(\down{b})$.
Since $A$ is finite, we can take the join $a=\bigvee f^{-1}(\down{b})$ of all elements of $f^{-1}(\down{b})$, and
from the above considerations it follows that $f^{-1}(\down{b})$ is the principal ideal
generated by $a$.
\end{proof}

In the next theorem we give a canonical bijection between the sets
$\Homj(A,B)$ and $\Homm^{\bot,1}(B_\bot,A_\bot)$, which will be the main tool
for the promised characterization of the operations in the centralizer
of a finite join\hyp{}semilattice.
Recall that $B_\bot$ and $A_\bot$ are lattices, and $\Homm^{\bot,1}(B_\bot,A_\bot)$
denotes the set of all meet\hyp{}homomorphisms
$g\colon B_\bot \to A_\bot$ satisfying $g(\bot)=\bot$ and $g(1)=1$.

\begin{theorem} \label{thm join-hom vs meet-hom}
Let $A, B$ be finite join\hyp{}semilattices, and
for every $f\in \Homj(A,B)$ and 
$g\in \Homm^{\bot,1}(B_\bot,A_\bot)$, let us define the maps 
$f^\triangleleft\colon B_\bot \to A_\bot$ and
$g^\triangleright\colon A \to B$ by
\[
f^\triangleleft(b) = 
  \begin{cases} 
   \bigvee f^{-1}(\down{b}), & \text{if } f^{-1}(\down{b}) \neq  \emptyset, \\
   \bot, & \text{if } f^{-1}(\down{b}) = \emptyset,
  \end{cases}
\qquad\quad
g^\triangleright(a) = \bigwedge g^{-1}(\up{a}).
\]
Then the following two maps are mutually inverse bijections:
\begin{align*}
\Homj(A,B) \to \Homm^{\bot,1}(B_\bot,A_\bot),\quad
&f \mapsto f^\triangleleft,\\
\Homm^{\bot,1}(B_\bot,A_\bot) \to \Homj(A,B),\quad
&g \mapsto g^\triangleright.
\end{align*}
\end{theorem}
\begin{proof}
First let us show that if $f$ is a join\hyp{}homomorphism from $A$ to $B$, then
$f^\triangleleft$ is a meet\hyp{}homomorphism from $B_\bot$ to $A_\bot$.
By Lemma~\ref{lemma inverse image of ideal}, $f^{-1}(\down{b})$ is either
empty or a principal ideal, and in the latter case $f^\triangleleft(b)$ is the greatest element of $f^{-1}(\down{b})$
by the definition of $f^\triangleleft$.
Therefore, we can reformulate the definition of  
$f^\triangleleft$ as follows:
\begin{equation} \label{eq f(a)<b <==> a<ft(b)}
\forall a \in A \; \forall b \in B_\bot \colon \  f(a) \leq b \iff a \leq f^\triangleleft(b).
\end{equation}
(Note that if $f^{-1}(\down{b}) = \emptyset$, then $f(a) \leq b$ 
does not hold for any $a \in A$. In this case we have $f^\triangleleft(b) = \bot$, and  
the only element $a \in A_\bot$ satisfying 
$a \leq f^\triangleleft(b) = \bot$ is $a = \bot$, thus the right hand side of 
\eqref{eq f(a)<b <==> a<ft(b)} does not hold for any $a \in A$ either.)
From \eqref{eq f(a)<b <==> a<ft(b)} we can deduce the following chain of 
equivalences for all $a \in A$ and $b_1,b_2 \in B_\bot$:
\begin{align*}
a \leq f^\triangleleft(b_1 \land b_2) 
&\iff f(a) \leq b_1 \land b_2\\
&\iff f(a) \leq b_1 \text{ and } f(a) \leq b_2\\
&\iff a \leq f^\triangleleft(b_1) \text{ and } a \leq f^\triangleleft(b_2)\\
&\iff a \leq f^\triangleleft(b_1) \land f^\triangleleft(b_2).
\end{align*}
Thus we have $a \leq f^\triangleleft(b_1 \land b_2) \iff a \leq f^\triangleleft(b_1) \land f^\triangleleft(b_2)$
for every $a \in A$, and this implies that 
$f^\triangleleft(b_1 \land b_2) = f^\triangleleft(b_1) \land f^\triangleleft(b_2)$, i.e., $f^\triangleleft$ is a meet\hyp{}homomorphism. 
To verify $f^\triangleleft(1_B)=1_A$, we just need to observe that 
$f^{-1}(\down{1_B})=f^{-1}(B)=A$, and the greatest element of $A$ is 
indeed $1_A$. Since $f^{-1}(\down{\bot})=\emptyset$, we have $f^\triangleleft(\bot)=\bot$.
This completes the proof of the claim $f^\triangleleft\in \Homm^{\bot,1}(B_\bot,A_\bot)$.

Now assume that $g \colon B_\bot \to A_\bot$ is a meet\hyp{}homomorphism such that
$g(1_B)=1_A$ and $g(\bot)=\bot$. Then $1_B \in g^{-1}(\up{a})$ for all $a \in A$, 
hence $g^{-1}(\up{a})$ is never empty.
Therefore, the dual of Lemma~\ref{lemma inverse image of ideal} shows that
$g^{-1}(\up{a})$ is a principal filter in $B_\bot$;
moreover, $g(\bot)=\bot$ implies that $\bot \notin g^{-1}(\up{a})$ 
for every $a \in A$. This shows that the map 
$g^\triangleright\colon A \to B, a \mapsto \bigwedge g^{-1}(\up{a})$ 
is well defined. One can prove, by an argument similar to that of the 
previous paragraph, that $g^\triangleright$ is a join\hyp{}homomorphism from $A$ to $B$.

It remains to prove that the maps
$f \mapsto f^\triangleleft$ and $g \mapsto g^\triangleright$
are inverses of each other. This follows immediately from the fact that
$g=f^\triangleleft$ and $f=g^\triangleright$ are both equivalent to 
\begin{equation} \label{eq f(a)<b <==> a<g(b)}
\forall a \in A \; \forall b \in B_\bot \colon \  f(a) \leq b \iff a \leq g(b).
\end{equation}
for all $f \in \Homj(A,B)$
and $g \in \Homm^{\bot,1}(B_\bot,A_\bot)$.
(Let us mention that a pair $(f,g)$ of maps satisfying 
\eqref{eq f(a)<b <==> a<g(b)} is called a 
\emph{monotone Galois connection}.)
\end{proof}

\begin{remark}
Let us give a categorical interpretation of 
Theorem~\ref{thm join-hom vs meet-hom}.
Let $\mathcal{J}$ denote the category of finite join-semilattices 
(with join-homomorphisms), and let $\mathcal{L}$ denote the category 
of finite lattices (with meet-homomorphisms). Then the following two
maps are mutually inverse functors, thus $\mathcal{J}$ and 
$\mathcal{L}$ are isomorphic categories:
\begin{align*}
F&\colon \mathcal{J} \to \mathcal{L},\; A \mapsto A_\bot,\; f \mapsto f^\triangleleft;\\ 
G&\colon \mathcal{L} \to \mathcal{J},\; B \mapsto B\setminus\{ 0_B \},\; g \mapsto g^\triangleright.
\end{align*}
(Here $B\setminus\{ 0_B \}$ is the join-semilattice obtained by 
removing the bottom element of the lattice $B$.
Of course, if $B$ is given as $B=F(A)=A_\bot$, then $0_B=\bot$.)
\end{remark}

Theorem~\ref{thm join-hom vs meet-hom} can be useful if $A$ is (much) larger
than $B$, as in this case it might be an easier task to determine the meet\hyp{}homomorphisms 
from $B_\bot$ to $A_\bot$ than describing the join\hyp{}homomorphisms from $A$ to $B$. 
This is the case when $A=S^n$ and $B=S$, where $S$ is a finite 
join\hyp{}semilattice: as mentioned before, the $n$-ary operations in 
$[\lor]^\ast$ are the join\hyp{}homomorphisms from $S^n$ to $S$, and these 
can be described in terms of the $1$- and $\bot$-preserving 
meet\hyp{}homomorphisms from $S_\bot$ to $(S^n)_\bot$, with the help of 
Theorem~\ref{thm join-hom vs meet-hom}. 
We formulate this characterization in the next corollary, and we complement it with 
the necessary and sufficient condition for the operation to depend on 
all of its variables.

\begin{corollary} \label{cor centralizer of semilattice}
Let $S=(S; \lor)$ be a finite semilattice, and let $n$ be a natural number.
The $n$-ary members of $[\lor]^\ast$ are exactly the operations $f$ of the form 
\[
f\colon S^n \to S, \  \x \mapsto \textstyle\bigwedge g^{-1}(\up{\x}),
\]
where $g\in \Homm^{\bot,1}(S_\bot,(S^n)_\bot)$; here $g$ is uniquely determined by $f$.
The operation $f$ depends on all of its variables if and only if
for each $i \in \{1, \dots, n\}$, the range of $g$ contains an element of $S^n$ whose $i$-th component is different from $1$.
If $S$ has a least element (i.e., if $S$ is a lattice), then the latter condition is satisfied
if and only if the range of $g$ contains a tuple from $(S\setminus\{1\})^n$.
\end{corollary}

\begin{proof}
An $n$-ary operation $f \in \O_S$ belongs to $[\lor]^\ast$ if and only if $f$
is a join\hyp{}homomorphism from $S^n$ to $S$. 
Applying Theorem~\ref{thm join-hom vs meet-hom} with $A=S^n$ and $B=S$, we see that these 
operations can be uniquely written as 
$f(\x)=g^\triangleright(\x)=\bigwedge g^{-1}(\up{\x})$ with
$g \in \Homm^{\bot,1}(B_\bot,A_\bot)$.

We prove that $f$ depends on its $i$-th variable if and only if the range
of $g$ contains a tuple $\s_i \in S^n$ such that the $i$-th component of
$\mathbf{s}_i$ is not equal to $1$. 
First assume that $f$ depends on the $i$-th variable; this means that 
there exist elements $\a,\a' \in S^n$ differing only in their $i$-th component
such that $b:=f(\a)$ and  $b':=f(\a')$ are different. 
We can assume without loss of generality that either $b < b'$ 
or $b$ and $b'$ are incomparable. In both cases we can conclude that
$\a \in f^{-1}(\down{b})$ and 
$\a' \notin f^{-1}(\down{b})$. 
By Theorem~\ref{thm join-hom vs meet-hom}, we have $g=f^\triangleleft$, 
thus $g(b)=\bigvee f^{-1}(\down{b})$. Therefore, 
$\a \leq g(b) $ and $\a' \nleq g(b)$.
This implies that the $i$-th component of the tuple $\s_i:=g(b) \in S^n$
(which certainly belongs to the range of $g$) is strictly less than $1$.

Conversely, let us suppose that there is an element $b \in S$ such that
the $i$-th component of $\s_i:=g(b)$ is less than $1$. Letting
$\s'_i$ be the tuple obtained from $\s_i$ by changing its
$i$-th component to $1$, we have $\s_i<\s'_i$.
Now $b \in g^{-1}(\up{\s_i})$ but $b \notin g^{-1}(\up{\s'_i})$,
therefore $g^{-1}(\up{\s_i}) \neq g^{-1}(\up{\s'_i})$.
Since $f=g^\triangleright$, this implies that
\[
f(\s_i)=g^\triangleright(\s_i)=\bigwedge g^{-1}(\up{\s_i})
\neq
\bigwedge g^{-1}(\up{\s'_i})=g^\triangleright(\s'_i)=f(\s'_i).
\]
Taking into account that $\s_i$ and $\s'_i$ differ only at the $i$-th component, 
we can conclude that $f$ does depend on its $i$-th variable.

If $S$ is a lattice, then $\mathbf{s}_1\land\dots\land\mathbf{s}_n$ is a tuple
in the range of $g$, and all of its components are less than $1$.
\end{proof}

\subsection{Counting}
Using the characterizations presented in the previous subsection, we can determine the exact number of $n$-ary operations in the centralizers of certain semilattices.
First we count the essentially $n$-ary operations commuting with the join operation of the smallest non-lattice semilattice.

\begin{proposition} \label{prop V}
Let $S=(\{a,b,1\},\lor)$ be the join semilattice with $a \lor b=1$. 
The number of essentially $n$-ary operations 
in the centralizer of $S$ is $8^n-6^n+2\cdot2^n+0^n$.
\end{proposition}

\begin{proof}
By Corollary~\ref{cor centralizer of semilattice}, we need to count the
elements $g\in\Homm^{\bot,1}(S_\bot,(S^n)_\bot)$ such that the
range of $g$ contains a tuple from $S^{i-1}\times\{a,b\}\times S^{n-i}$
for every $i\in\{1,\ldots,n\}$.
For an arbitrary map $g\colon S_\bot \to (S^n)_\bot$, we have
$g\in\Homm^{\bot,1}(S_\bot,(S^n)_\bot)$ if and only if $g(1)=1$, $g(\bot)=\bot$
and $g(a) \land g(b) = \bot$; moreover, such a map is uniquely determined by $g(a)$ and $g(b)$.
We distinguish four cases upon these values 
(we denote by $f=g^\triangleright$ the element of $[\lor]^\ast$ corresponding to $g$ in
Corollary~\ref{cor centralizer of semilattice}).

\begin{enumerate}[label=(\arabic*)]
\item If $g(a)=\bot=g(b)$, then $f$ is constant $1$, hence $f$ is essentially $n$-ary
if and only if $n=0$. Thus the number of essentially $n$-ary operations of this
type is $0^n$, i.e., it is $1$ if $n=0$  and $0$ if $n>0$ 
(here it is convenient to use the convention $0^0=1$; for more justification, see \cite{Knuth}).

\item If $g(a)\neq\bot=g(b)$, then $f$ depends on all of its variables if and only if
$g(a) \in \{a,b\}^n$, thus the number of essentially $n$-ary operations 
$f \in [\lor]^\ast$ of this type is $2^n$ (for $n=0$ we get the constant $a$ function). 

\item If $g(a)=\bot\neq g(b)$, then, similarly to the previous case,
we have $2^n$ functions (here $n=0$ corresponds to the constant $b$ function).

\item If $g(a)\neq\bot\neq g(b)$, then let $\s:=g(a)\in S^n$ and $\mathbf{t}:=g(b)\in S^n$.
Writing these two tuples below each other, we get the $2 \times n$ matrix
$
\bigl(
\begin{smallmatrix}
s_1 & \dots & s_n \\
t_1 & \dots & t_n
\end{smallmatrix}
\bigr).
$
Now $f$ depends on all of its variables if and only if no column of this matrix
is $(1,1)^\mathsf{T}$, and there are $8^n$ such matrices. However, some of
the corresponding maps $g$ will violate the condition $g(a) \land g(b) = \bot$:
we must exclude those matrices that contain neither $(a,b)^\mathsf{T}$ nor 
$(b,a)^\mathsf{T}$ as a column. The number of such matrices is $6^n$, so
we obtain $8^n-6^n$ essentially $n$-ary operations $f \in [\lor]^\ast$ in this case.
\end{enumerate}
Summing up the four cases, we see that the number of essentially $n$-ary operations 
in $[\lor]^\ast$ is $8^n-6^n+2\cdot2^n+0^n$.
\end{proof}

\begin{remark} \label{remark all n-ary}
It is easy to see that if the number of essentially $n$-ary 
operations in a clone $C$ is $p_n$, then the number of all 
operations of arity $n$ in $C$ is 
$\sum_{k=0}^{n} \binom{n}{k}p_k$.
Thus, by Proposition~\ref{prop V} (and by the binomial theorem), 
the number of $n$-ary operations in the centralizer of the 
join operation of the semilattice $(\{a,b,1\},\lor)$ is
\[
\sum_{k=0}^{n} \mbinom{n}{k} (8^k-6^k+2\cdot2^k+0^k) =
9^n-7^n+2\cdot3^n+1^n.
\]
\end{remark}

Next we provide a general formula for the number of essentially $n$-ary 
operations commuting with the join operation of a finite lattice, and 
then we apply it to the case of finite chains.

\begin{theorem} \label{thm size of centralizer of semilattice}
Let $S=(S; \lor, \land)$ be a finite lattice, and let $n$ be a natural number. 
The number of essentially $n$-ary operations in $[\lor]^\ast$ is
\[
\sum_{b \in S} (\lvert \Homj^{0}(S,\up{b})\rvert-1)^n =
\sum_{b \in S} (\lvert \Homm^{1}(\up{b},S)\rvert-1)^n.
\]
\end{theorem}

\begin{proof}
According to Theorem~\ref{thm Larose}, the essentially $n$-ary members
of the centralizer are in a one-to-one correspondence with the tuples 
$(u_1, \dots, u_n) \in \Homj(S,S)^n$ such that $u_1(0) = \dots = u_n(0)$
and none of the $u_i$ are constant. Let $b:=u_1(0)$, then each $u_i$
maps $S$ to $\up{b}$ in such a way that the least element of $S$ is mapped to
the least element of $\up{b}$, i.e., 
$u_i \in \Homj^{0}(S,\up{b})$ for $i=1,\dots,n$.
However, we need to exclude the constant $b$ function, hence the number of
choices for each $u_i$ is $\lvert \Homj^{0}(S,\up{b})\rvert-1$,
and this gives the first formula. (Note that for $b=1$, the principal filter
$\up{b}$ has just one element, thus $\lvert \Homj^{0}(S,\up{b})\rvert - 1 = 0$.
Therefore, the contribution of $b=1$ to the sum is $0^n$, and this could
be omitted if $n>0$. However, for $n=0$, we need to keep the term $0^0=1$
in order to get the correct number of nullary operations,
which is clearly $|S|$.)

The second formula follows from the first one by applying 
Theorem~\ref{thm join-hom vs meet-hom} to $A=S$ and $B=\up{b}$.
If $f \in \Homj(S,\up{b})$ and $g \in \Homm^{\bot,1}((\up{b})_\bot,S_\bot)$
correspond to each other under the bijections of 
Theorem~\ref{thm join-hom vs meet-hom}, then $f(0)=b$ if and only if
$g(0) \neq \bot$ (i.e., $g$ does not take the value $\bot$ except for
$g(\bot)=\bot$). Therefore, we obtain a bijection from
$\Homj^{0}(S,\up{b})$ to $\Homm^{1}(\up{b},S)$ by restricting 
$f^\triangleleft$ to the set $\up{b}$ for each $f \in \Homj^{0}(S,\up{b})$.
\end{proof}

\begin{remark}
The second formula of the above theorem can be also derived 
directly from Corollary~\ref{cor centralizer of semilattice} as follows.
We need to count the meet\hyp{}homomorphisms $g \in \Homm^{\bot,1}(S_\bot,(S^n)_\bot)$
whose range satisfies the conditions of 
Corollary~\ref{cor centralizer of semilattice}.
If $S$ is a lattice and $g$ is such a meet\hyp{}homomorphism, then 
there is a least element $b \in S$ such that $g(b)\neq\bot$.
Restricting $g$ to $\up{b}$, we get a $1$-preserving meet\hyp{}homomorphism 
from $\up{b}$ to $S^n$, which can be viewed as an $n$-tuple $(g_1,\dots,g_n)$ of 
$1$-preserving meet\hyp{}homomorphisms from $\up{b}$ to $S$.
The range of $g$ contains a tuple from $(S \setminus \{ 1 \})^n$
if and only if $g(b) \in (S \setminus \{ 1 \})^n$, which holds if and only if 
none of the $g_i$ are constant $1$.
Therefore, there are exactly $(\lvert \Homm^{1}(\up{b},S)\rvert-1)^n$ 
such tuples $(g_1,\dots,g_n)$, and this proves the second formula of
Theorem~\ref{thm size of centralizer of semilattice}.
This argument does not work if $S$ is not a lattice, even though
Corollary~\ref{cor centralizer of semilattice} holds in that case, too.
The problem is that there may be no least element $b$ with $g(b)\neq\bot$;
in fact, $g^{-1}(S^n)$ is not necessarily closed under meets.
Nevertheless, as we have seen in Proposition~\ref{prop V}, 
Corollary~\ref{cor centralizer of semilattice} can be used to count
the essentially $n$-ary operations in a semilattice even if it is not a lattice.
\end{remark}

\begin{proposition} \label{prop chain}
The number of essentially $n$-ary operations commuting with 
the join\hyp{}operation of a chain of cardinality $\ell$ is
\[
\sum_{i=1}^\ell \Bigl[\mbinom{\ell+i-2}{\ell-1}-1\Bigr]^n.
\]
\end{proposition}

\begin{proof}
First, as an auxiliary result, let us count the join\hyp{}homomorphisms
from an $r$-element chain $A=\{a_1<\dots<a_r\}$ to an 
$s$-element chain $B=\{b_1<\dots<b_s\}$. Clearly, the
join\hyp{}homomorphisms in this case are just the monontone maps, 
thus an element of $\Homj(A,B)$ can be given
by a nondecreasing sequence $f(a_1) \leq\dots\leq f(a_r)$ in $B$.
These sequences can be viewed as $r$-combinations with repetitions
from the elements $b_1,\dots,b_s$, and the number of such 
combinations is $\binom{s+r-1}{r}=\binom{s+r-1}{s-1}$.

Now let $S$ be a chain of cardinality $\ell$, and let $b \in S$.
The $0$-preserving join\hyp{}homomorphisms  from $S$ to $\up{b}$ 
are in a one-to-one correspondence with the join\hyp{}homomorphisms 
from $S\setminus\{ 0 \}$ to $\up{b}$ 
(by restricting to $S\setminus\{ 0 \}$).
Note that $S\setminus\{0\}$ is a chain of size $\ell-1$, 
and if $b$ is the $i$-th element from the top in $S$, then
$\up{b}$ is an $i$-element chain. Thus, by the considerations
made in the first paragraph, we have 
\[
\lvert \Homj^{0}(S,\up{b})\rvert = 
\lvert \Homj(S\setminus\{0\},\up{b})\rvert =
\mbinom{\ell+i-2}{\ell-1}.
\]
Applying Theorem~\ref{thm size of centralizer of semilattice}
completes the proof: we just need to substitute the above formula
into the first sum of Theorem~\ref{thm size of centralizer of semilattice}
(replacing the summation variable $b$ by $i$).
\end{proof}

\begin{remark}
For an arbitrary poset $A$, the number of monotone maps
from $A$ to an $s$-element chain is a polynomial in $s$, called
the \emph{order polynomial} of $A$. As we have seen in the
proof of the above proposition, the order polynomial of the 
$r$-element chain is $\binom{s+r-1}{r}$. 
This fact, and much more about order polynomials can be
found in \cite{Stanley}.
\end{remark}

\begin{remark} \label{remark all n-ary for chains}
Similarly to Remark~\ref{remark all n-ary}, we can
derive from Proposition~\ref{prop chain} that
the number of $n$-ary operations in the centralizer of the 
join (or meet) operation of an $\ell$ element chain is
\[
\sum_{k=0}^{n} \mbinom{n}{k} \sum_{i=1}^\ell \Bigl[\mbinom{\ell+i-2}{\ell-1}-1\Bigr]^k =
\sum_{i=1}^\ell \sum_{k=0}^{n}  \mbinom{n}{k}\Bigl[\mbinom{\ell+i-2}{\ell-1}-1\Bigr]^k =
\sum_{i=1}^\ell \mbinom{\ell+i-2}{\ell-1}^n.
\]
\end{remark}

\begin{example}
For the two-element chain (regarded as a semilattice), 
Proposition~\ref{prop chain} gives the formula $1^n+0^n$ 
for the number of essentially $n$-ary operaions (recall that $0^0$ is defined as $1$),
and Remark~\ref{remark all n-ary for chains} gives $2^n+1^n$ 
for the number of all $n$-ary operations in the centralizer
(which are of course also easily verified without our results).
Similarly, for the three-element chain, we have
$5^n+2^n+0^n$ essentially $n$-ary operations and
$6^n+3^n+1^n$ operations of arity $n$ in the centralizer of the join operation;
for the four-element chain we get the numbers
$19^n+9^n+3^n+0^n$ and $20^n+10^n+4^n+1^n$, etc.
A formula for the number of $n$-ary operations in the centralizer of the meet operation of the three-element chain
appeared already in \cite{MachidaRosenberg}:
\[
3^{n+1} + \sum_{\substack{0 \leq p < n \\ 0 \leq q \leq n-p}} \mbinom{n}{p} \mbinom{n-p}{q} (3^p 2^q -1),
\]
This can be simplified to $6^n+3^n+1^n$ with the help of
the binomial theorem.
\end{example}

\newpage

\section{Centralizers of finite distributive lattices} \label{sect:distributive}

In this section $L=(L; \lor, \land)$ denotes a finite distributive lattice with greatest element $1$ and least element $0$.
We use the symbol $\2$ for the two-element lattice $\{0,1\}$, and $\mathcal{P}(U)$ denotes the lattice of all subsets of a set $U$.
Note that $\mathcal{P}(\{1, \dots, n\})\cong\2^n$.
\begin{theorem} \label{thm f=u_1*...*u_n lattice}
 Let $L=(L; \lor, \land)$ be a finite distributive lattice and $f \in \OL^{(n)}$. Then the following are equivalent:
	\begin{enumerate}
		\item $f \in [\lor, \land]^\ast$;
		\item there exist  unary operations $u_1, \dots, u_n \in [\lor, \land]^\ast$ such that $f(x_1,  \dots, x_n) = u_1(x_1) \lor \dots \lor u_n(x_n)$ and for all $i, j \in \{1, \dots, n\}, i\neq j$ we have $u_i(1) \land u_j(1) = u_1(0)= \dots = u_n(0)$.
	\end{enumerate}
\noindent Furthermore, the operation $f$ given by $(2)$ depends on all of its variables if and only if none of the unary operations $u_i$ are constant.
\end{theorem}
\begin{proof}
(1) $\Rightarrow$ (2): 
 As in the proof of Theorem~\ref{thm Larose}, we define the unary operations $u_1, \dots, u_n$ as $u_1(x)=f(x, 0, \dots, 0), \dots, u_n(x)=f(0, \dots, 0, x)$. By the theorem, for these unary operations we have $f(x_1, \dots, x_n) = u_1(x_1) \lor \dots \lor u_n(x_n)$ and $u_1(0)= \dots = u_n(0) = f(0, \dots, 0)$. Since $f$ and the constant $0$ operation belong to $[\lor, \land]^\ast$, we have $u_1, \dots, u_n \in [\lor, \land]^\ast$. It only remains to show that for all $i, j \in \{1, \dots, n\}, i\neq j$ we have $u_i(1) \land u_j(1) = f(0, \dots, 0)$. For notational simplicity, let us assume that $i=1$ and $j=2$; the proof of the general case is similar. Using that $f$ commutes with the operation $\land$, with the help of the $2$ by $n$ matrix
$\bigl( 
\begin{smallmatrix}
1 & 0 & 0 & \ldots & 0 \\
0 & 1 & 0 & \ldots & 0
\end{smallmatrix}
 \bigr)$
\noindent we can conclude that the following equality holds:
\begin{align*} f(0, \dots, 0) &= f\big(1 \land 0,\ 0 \land 1,\ 0 \land 0,\ \dots,\ 0 \land 0) \\
&= f(1, 0, 0, \dots, 0) \land f(0, 1, 0, \dots, 0) = u_1(1) \land u_2(1).
\end{align*}

(2) $\Rightarrow$ (1): Since $\lor \in [\lor]^\ast$ and $u_1, \dots, u_n \in [\lor, \land]^\ast$, we have $f(x_1, \dots, x_n)=u_1(x_1) \lor \dots \lor u_n(x_n) \in [\lor]^\ast$. Therefore, to complete the proof we have to show that $f$ commutes with $\land$: 
\[
f(x_1, \dots, x_n) \land f(y_1, \dots, y_n)=f( x_1 \land y_1, \dots, x_n \land y_n).
\]
Thus we need to prove that for all $x_1,\dots,x_n,y_1,\dots,y_n$ we have 
\[
\big( u_1(x_1) \lor \dots \lor u_n(x_n) \big) \land \big(u_1(y_1) \lor \dots \lor u_n(y_n) \big) = u_1(x_1 \land y_1) \lor \dots \lor u_n(x_n \land y_n).
\]
Using the notation $c_i:=u_i(x_i), d_i:=u_i(y_i) (i=1,\dots,n)$, the above equality can be written as
\[
( c_1 \lor \dots \lor c_n ) \land ( d_1 \lor \dots \lor d_n ) =
u_1(x_1 \land y_1) \lor \dots \lor u_n(x_n \land y_n).
\]
Since $L$ is distributive, we have
\begin{align}\label{eq (c1+...+cn)(d1+...+dn)=+c_id_j}
( c_1 \lor \dots \lor c_n ) \land ( d_1 \lor \dots \lor d_n )
= \bigvee_{i,j =1}^{n} (c_i \land d_j).
\end{align}
From $u_1, \dots, u_n \in [\lor, \land]^\ast$ it follows that these operations are monotone, hence 
$u_1(0)=u_i(0) \leq c_i,d_i \leq u_i(1)$ for every $i$; moreover, (2) also implies that for all $i\neq j$, we have
\[
u_1(0) = u_i(0) \land u_j(0) \leq c_i\land d_j \leq u_i(1) \land u_j(1) =  u_1(0).
\]
Thus $c_i\land d_j=u_1(0) \leq c_i\land d_i$ whenever $i \neq j$, so we can omit $c_i\land d_j$ from the join on the right hand side of \eqref{eq (c1+...+cn)(d1+...+dn)=+c_id_j}, and using the fact that each $u_i$ commutes with $\land$ we obtain the desired equality:
\begin{align*}
( c_1 \lor \dots \lor c_n ) \land ( d_1 \lor \dots \lor d_n )
&= \bigvee_{i =1}^{n} (c_i \land d_i) \\
&= \bigvee_{i =1}^{n} (u_i(x_i) \land u_i(y_i)) \\
&= \bigvee_{i =1}^{n} (u_i(x_i \land y_i)).
\end{align*}
\end{proof}

\begin{lemma} \label{lemma essentially n-ary in C* <-> 2^n sublattice}
 Let $L=(L; \lor, \land)$ be a finite distributive lattice. Then the following are equivalent:
	\begin{itemize}[leftmargin=.6in]
	    	\item [(Ess)] there exists an essentially $n$-ary operation in  $[\lor, \land]^\ast$;
		\item [(Sub)] there exists a sublattice of $L$ that is isomorphic to $\2^n$.
	\end{itemize}
\end{lemma}
\begin{proof}
(Ess) $\Rightarrow$ (Sub): Let $f \in \OL$ be an essentially $n$-ary operation. Then by Theorem~\ref{thm f=u_1*...*u_n lattice} we have $f(x_1, \dots, x_n) = u_1(x_1) \lor \dots \lor u_n(x_n)$ and $u_i(1) \land u_j(1)=u_1(0)=\dots=u_n(0)$ for all $i \neq j$. Let us introduce the notation $a_i=u_i(1)$ and $b=u_1(0)$. Then for all $i \neq j$ we have $a_i \land a_j = b$ and since $f$ is essentially $n$-ary, we also have $a_i > b$.
From this it is not hard to
deduce using distributivity that $\mathcal{P}(\{1, \dots, n\})\hookrightarrow L,~I\mapsto
\bigvee\{a_{i}:i\in I\}$ is an embedding. (Alternatively, one can verify with
the help of Theorem~360 of \cite{Gr} that $\{a_{1},\dots,a_{n}\}$ is an
independent set in the sublattice $\up{b}$, hence it generates a sublattice isomorphic to $\mathcal{P}(\{1, \dots, n\})$.)

(Sub) $\Rightarrow$ (Ess): Let us suppose that there is a sublattice of $L$ isomorphic to $\2^n$, let $b$ be the least element and $a_i \, (i\in \{1, \dots, n\})$ be the atoms of this cube.  Then we have $b=a_i \land a_j$ for all $i \neq j$. Let us define the operations $u_1, \dots, u_n$ as $u_i(x_i):= (x_i \lor b) \land a_i$ for all $i \in \{1, \dots, n\}$. Since $L$ is distributive, $u_i$ is an endomorphism of $L$, i.e., we have $u_i \in [\lor, \land]^\ast$. Note that $u_i(0) = b$ and $u_i(1)=a_i$, therefore we have $u_1(0)= \dots = u_n(0)$ and also $u_i(1) \land u_j(1)=u_1(0)$ for all $i \neq j$. By Theorem~ \ref{thm f=u_1*...*u_n lattice} this means that the operation $f \in \OL$ defined as $f(x_1, \dots, x_n) := u_1(x_1) \lor \dots \lor u_n(x_n)$ belongs to $[\lor, \land]^\ast$. Since none of the $u_i$ are constant, Theorem~\ref{thm f=u_1*...*u_n lattice} also implies that $f$ is essentially $n$-ary.
\end{proof}	

\begin{corollary}
	For a finite distributive lattice $L=(L; \lor, \land)$ the following are equivalent:
	\begin{itemize}
	\item every operation in $[\lor, \land]^\ast$ is essentially at most unary;
	\item $L$ is a chain.
	\end{itemize}
\end{corollary}
	
	Although the next lemma follows from the description of projective and injective distributive lattices \cite{Balbes}, we provide a short proof.
	
\begin{lemma}
 Let $L=(L; \lor, \land)$ be a finite distributive lattice. Then the following are equivalent:
	\begin{itemize}[leftmargin=.6in]
	\item [(Sub)] there exists a sublattice of $L$ that is isomorphic to $\2^n$;
	\item [(Quo)] there exists a congruence $\teta$ of $L$  such that $L/ \teta$ is isomorphic to $\2^n$.
	\end{itemize}
\end{lemma}
\begin{proof}

Instead of $\2^n$, it will be more convenient to use the lattice $K_n:=\mathcal{P}(\{1,...,n\})$, which is clearly isomorphic to $\2^n$. To prove (Sub) $\Rightarrow$ (Quo), assume that $L$ has a sublattice that is isomorphic to $\2^n$. Identifying this sublattice with $K_n$, we
may assume without loss of generality that $K_n$ itself is a
sublattice of $L$. For any $i\in \{1, \dots, n\}$, the principal ideal generated by
$\{1, \dots, n\}\setminus\{i\}$ does not contain $\left\{  i\right\}  $, hence, by the
prime ideal theorem for distributive lattices, there is a prime ideal $P_i$
of $L$ that does not contain $\{i\}$ (see Corollary~116 in \cite{Gr}).
Consequently, there is a homomorphism $\varphi_{i}\colon L\rightarrow
\2$ mapping $P_{i}$ to $0$ and $L\setminus P_{i}$ to $1$. In particular, we have
$\varphi_{i}(\{i\})=1$ and $\varphi_{i}(\{j\})=0$ for all $j\neq i$. Combining
these maps we obtain a homomorphism $\varphi\colon L\rightarrow\2^{n},~x\mapsto\left(  \varphi_{1}\left(  x_{1}\right)  ,\ldots,\varphi
_{n}\left(  x_{n}\right)  \right)  $. We have $\varphi\left(  \{i\}\right)
=\left(  0,\ldots,0,1,0,\ldots,0\right)  $, where the $1$ appears in the
$i$-th coordinate. These elements generate $\2^{n}$, hence $\varphi$
is surjective, and this proves (Sub).

For (Quo) $\Rightarrow$ (Sub), let us suppose that $\teta$ is a congruence of $L$
such that $L/\teta$ is isomorphic to $K_n$. For every $I\in K_n$, let $C_{I}$ denote the congruence class of $\teta$ corresponding
to $I$ at this isomorphism. Let $a$ be the greatest element of $C_{\emptyset}%
$, and let $b_{i}$ be the least element of $C_{\{i\}}$ for all $i\in \{1, \dots, n\}$. Then $c_{i}:=a\vee b_{i}$ belongs to $C_{\{i\}}$, and $c_{i}\wedge c_{j}$
belongs to $C_{\emptyset}$ whenever $i\neq j$. Moreover, $c_{i}\wedge
c_{j}=(a\vee b_{i})\wedge(a\vee b_{j})\geq a$, hence $c_{i}\wedge c_{j}=a$, as
$a$ is the greatest element of $C_{\emptyset}$. From this it follows using the same argument as in the proof of Lemma~\ref{lemma essentially n-ary in C* <-> 2^n sublattice} that $K_n\rightarrow L,~I\mapsto
\bigvee\{c_{i}:i\in I\}$ is an embedding.
\end{proof}

The previous two lemmas together give us the following characterization of the existence of essentially $n$-ary operations in the centralizer of a finite distributive lattice.

\begin{theorem} \label{thm essentially n-ary in C* <-> 2^n sublattice <-> 2^n congruence}
 Let $L=(L; \lor, \land)$ be a finite distributive lattice. Then the following are equivalent:
	\begin{itemize}[leftmargin=.6in]
		\item [(Ess)] there exists an essentially $n$-ary operation in  $[\lor, \land]^\ast$;
		\item [(Sub)] there exists a sublattice of $L$ that is isomorphic to $\2^n$;
		\item [(Quo)] there exists a congruence $\teta$ of $L$  such that $L/ \teta$ is isomorphic to $\2^n$.
	\end{itemize}
\end{theorem}

\section{Centralizers of finite lattices}\label{sect:lattices}

 In this section our goal is to investigate whether the results proved in Section~\ref{sect:distributive} hold for arbitrary lattices.
 Let us look at Theorem~\ref{thm f=u_1*...*u_n lattice} first. Note that in the proof of $(1) \Rightarrow (2)$ we did not use that the lattice $L$ was distributive neither that $L$ was finite, and thus the proof provided  there shows that this implication holds for arbitrary bounded lattices. However, in the proof of $(2) \Rightarrow (1)$ we used distributivity, and the following example shows that this implication does not hold for arbitrary lattices.

 \begin{example}
 Let $L$ be the lattice shown on Figure~\ref{counterexample (2) -> (1)}. 
 \begin{figure}
    \centering
	\includegraphics[width=0.6\linewidth]{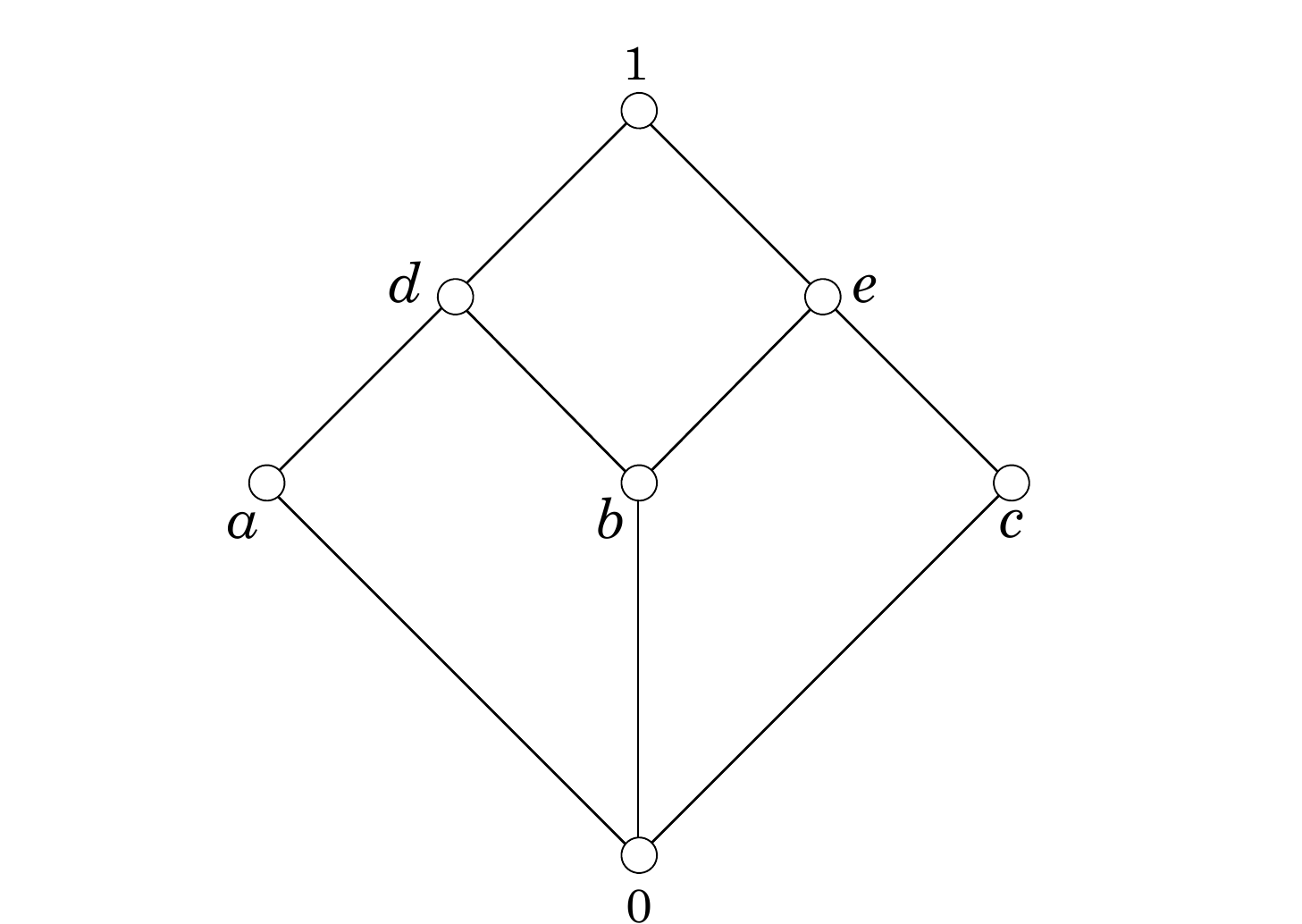}
	\caption{A nondistributive counterexample to the implication $(2) \Rightarrow (1)$ of Theorem~\ref{thm f=u_1*...*u_n lattice}.}  \label{counterexample (2) -> (1)}
\end{figure}
\noindent Let us define the operations $u_1$ and $u_2$ as
\[ u_1(x)=
\begin{cases}
0, & \text{ if } x \in \{0,b\}, \\
a, & \text{ if } x \in \{a, d\}, \\
b, & \text{ if } x \in \{c, e\}, \\
d, & \text{ if } x \in \{ 1 \};
\end{cases}
\quad \quad
u_2(x)=
\begin{cases}
0 & \text{ if } x \in \{ 0,a,b,d \}, \\
c & \text{ if } x \in \{c,e,1 \}.
\end{cases}\]

\noindent It is easy to check that $\ker(u_1)$ is a  congruence of $L$, and $u_1$ establishes an isomorphism from the quotient lattice $L / \ker(u_1) = \{ \{0, b\}, \{a,d\}, \{c,e\}, \{1\} \}$ to the sublattice $\{ 0, a, b, d\}$. Therefore $u_1$ is an endomorphism of $L$. Similarly, one can verify that $u_2 \in [\lor, \land]^\ast$. Let $f(x_1, x_2) = u_1(x_1) \lor u_2(x_2)$, we will show that $f$ does not belong to $[\land]^\ast$. Let us suppose that $f$ commutes with $\land$. Then applying the definition of commutation to the $2$ by $2$ matrix 
$\bigl( 
\begin{smallmatrix}
a & c \\
c & c
\end{smallmatrix}
 \bigr)$,
we have the following equality:
$$f(a \land c, c \land c) = f(a,c) \land f(c,c);$$
that is,
$$u_1(a \land c) \lor u_2(c \land c)=(u_1(a) \lor u_2(c)) \land (u_1(c) \lor u_2(c)).$$
However, the left hand side evaluates to $c$, and the value of the right hand side is $e$;
$$c = 0 \lor c=u_1(0) \lor u_2(c)= (a \lor c) \land (b \lor c)=1 \land e = e$$
which is a contradiction.
\end{example}

We have seen that the implication $(2) \Rightarrow (1)$ of Theorem~\ref{thm f=u_1*...*u_n lattice} does not hold for arbitrary lattices, but interestingly for finite simple lattices we have $(1) \Leftrightarrow (2)$. Moreover, as we shall see in the following remark, for a finite simple lattice $L=(L; \lor, \land)$, every operation in $[\lor, \land]^\ast$  depends on at most one variable.

\begin{remark} \label{remark (2)->(1) simple lattices}
If $L$ is a finite simple lattice, then $(2) \Rightarrow (1)$ holds even if $L$ is not distributive. Assume that $L$ is a finite simple lattice, $u_1, \dots ,u_n$ are endomorphisms of $L$ such that $u_1(0)= \dots =u_n(0)$ and $u_i(1)\land u_j(1)=u_1(0)$ whenever $i\neq j$, and let $f(x_1, \dots ,x_n)=u_1(x_1)\lor \dots \lor u_n(x_n)$. Simplicity of $L$ implies that the kernel of every $u_i$ is one of the two relations $L^2$ or $\{(a,a) \mid a \in L\}$; therefore, every $u_i$ is either a constant operation or an automorphism of $L$. 
We distinguish three cases on the number of automorphisms occurring in $u_1(x_1)\lor\dots \lor u_n(x_n)$.
\begin{enumerate}[label=(\roman*)]
\item If each $u_i$ is constant, then $f$ is also constant, hence $f \in [\lor,\land]^\ast$.

\item If $u_i$ and $u_j$ are automorphisms with $i\neq j$, then $1=u_i(1) \land u_j(1) = u_i(0)=0$, which is a contradiction. 

\item In the remaining cases we can assume without loss of generality that $u_1$ is an automorphism and $u_2, \dots ,u_n$ are all constants. Then $u_1(0)=0$, thus $u_2(0)= \dots =u_n(0)=0$ and $f(x_1, \dots ,x_n)=u_1(x_1)$. 
\end{enumerate}
\noindent Therefore, if $(2)$ holds for a finite simple lattice, then $f$ depends on at most one variable, and $f$ is equivalent to an automorphism or a constant, hence $f \in [\lor, \land]^\ast$.
Thus $(1) \Leftrightarrow (2)$ holds for finite simple lattices; moreover, every operation in $[\lor,\land]^\ast$ depends on at most one variable.
\end{remark}

Before investigating further which results of Section~\ref{sect:distributive} hold for arbitrary finite lattices, we need to recall some facts from universal algebra. First we define the so-called product congruence; this definition can also be found in \cite{BS}.

\begin{definition}
Let $A_1, \dots, A_n$ be algebras of the same type and let $\teta_i \in \Con(A_i)$ for all $i \in \{1, \dots, n\}$. The \emph{product congruence}
$\teta = \teta_1 \times \dots \times \teta_n$
on $A_1 \times \dots \times A_n$ is defined by
$\big((a_1, \dots, a_n) , (b_1, \dots, b_n)\big) \in \teta \Leftrightarrow \forall i \in \{1, \dots, n\}\colon (a_i, b_i) \in \teta_i$.
\end{definition}

Next we give two lemmas about congruences of direct products in congruence distributive varieties. The first one is a special case of Lemma~11.10 of \cite{BS}; the second one is implicit in \cite{CM,CGL}, but we include the proof for the sake of self-containedness. 
The variety of lattices is congruence distributive, hence we can use them in our study of centralizers of lattices.  
These two lemmas will also be helpful later in describing the centralizers of the clones over the two-elment set (see Section~\ref{sect:boole}).

\begin{lemma}[\cite{BS}] \label{lemma CD variety theta=theta_1 x ... x theta_n}
Let $\mathcal{V}$ be a congruence distributive variety, $A_1, \dots, A_n \in \mathcal{V}$ and $\teta \in \Con(A_1 \times \dots \times A_n)$. Then there exist $\teta_i \in \Con(A_i)$ for all $i \in \{1, \dots, n\}$ such that $\teta = \teta_1 \times \dots \times \teta_n$.
\end{lemma}

\begin{lemma} \label{lemma CD variety essentially n-ary in C* <-> A_1/theta_1 x ... x A_n/theta_n <= A}
Let $\mathcal{V}$ be a congruence distributive variety and $A \in \mathcal{V}$. Then the following are equivalent:
\begin{enumerate}[label=(\roman*)]
	\item There exists an essentially $n$-ary operation $f \in \OA$ that is a homomorphism from $A^n$ to $A$.
	\item There exist $\teta_i \in \Con(A)$ such that $\teta_i\neq A^2$ for all $i \in \{1, \dots, n\}$, and $A_1/\teta_1 \times \dots \times A_n/\teta_n$ embeds into $A$.
\end{enumerate}
\end{lemma}
\begin{proof}
(i) $\Rightarrow$ (ii): Let $f \colon A^n \rightarrow A$ be a homomorphism. Then by Lemma~\ref{lemma CD variety theta=theta_1 x ... x theta_n} there exist $\teta_1, \dots, \teta_n \in \Con(A)$ such that $\ker(f)=\teta_1 \times \dots \times \teta_n$. Since $f$ is essentially $n$-ary, we have that $\teta_i \neq A^2$ for all $i \in \{1, \dots, n\}$. By the homomorphism theorem we have $A^n/\ker(f)=A^n / (\teta_1 \times \dots \times \teta_n) = A/\teta_1 \times \dots \times A/\teta_n \cong f(A^n)$, and since $f(A^n)$ is a subalgebra of $A$, it is clear that $A_1/\teta_1 \times \dots \times A_n/\teta_n$ is embeddable into $A$.

(ii) $\Rightarrow$ (i): Let $\phi \colon A/\teta_1 \times \dots \times A/\teta_n \hookrightarrow A$ be an embedding, let $\teta \in \Con(A^n)$ be the product congruence $\teta = \teta_1 \times \dots \times \teta_n$  and  let $\nu$ denote the natural homomorphism from $A^n$ to $A^n/\ker\phi = A^n/\teta$. Then we define the operation $f \in \OA^{(n)}$ as $f(x_1, \dots, x_n)=\phi(x_1/\teta_1, \dots , x_n / \teta_n)=\phi (\nu (x_1, \dots, x_n))$. Thus $f=\phi \circ \nu$ is a homomorphism, and since $\teta_i \neq A^2$ for all $i$, we have that $f$ is essentially $n$-ary.
\end{proof}

The following corollary of Lemma~\ref{lemma CD variety essentially n-ary in C* <-> A_1/theta_1 x ... x A_n/theta_n <= A} gives an upper bound for the essential arity of operations in the centralizer of a finite algebra in a congruence distributive variety.

\begin{corollary}\label{cor:log2}
 Let $\mathcal{V}$ be a congruence distributive variety, let $A \in \mathcal{V}$ be a finite algebra and let $C$ denote the clone of term operations of $A$. If there is an essentially $n$-ary operation in $C^\ast$, then we have $n \leq \log_2 |A|$. In other words, the essential operations in $C^\ast$ are at most $\log_2 |A|$-ary.
\end{corollary}
\begin{proof}
By Lemma~\ref{lemma CD variety essentially n-ary in C* <-> A_1/theta_1 x ... x A_n/theta_n <= A}, if we have an essentially $n$-ary operation in $C^\ast$, then there exist $\teta_1, \dots, \teta_n \in \Con(A)$ such that $A/\teta_1 \times \dots \times A/\teta_n$ is embeddable into $A$. Therefore we have $| A/\teta_1 \times \dots \times A/\teta_n | = | A/\teta_1 | \cdot \ldots \cdot| A/\teta_n |\leq |A|$, and since for every $i$ we have $\teta_i \neq A^2$, it follows that $2^n=2 \cdot \ldots \cdot 2 \leq | A/\teta_1 | \cdot \ldots \cdot| A/\teta_n |\leq |A|$. Thus $n=\log_2(2^n) \leq \log_2(|A|)$.
\end{proof}

Now we will focus on Theorem~\ref{thm essentially n-ary in C* <-> 2^n sublattice <-> 2^n congruence}, or more precisely, we investigate which implications between (Ess), (Sub) and (Quo) hold for arbitrary finite lattices.
Using that the variety of lattices is congruence distributive, and also that every lattice (except the one-element lattice) has a two-element sublattice, as a corollary of Lemma~\ref{lemma CD variety essentially n-ary in C* <-> A_1/theta_1 x ... x A_n/theta_n <= A} we can conclude that the implication (Ess) $\Rightarrow$ (Sub) of Theorem~\ref{thm essentially n-ary in C* <-> 2^n sublattice <-> 2^n congruence} also holds for arbitrary finite lattices.

However, (Sub) $\Rightarrow$ (Ess) does not hold in general; partition lattices provide counterexamples to this implication. Indeed, every finite lattice (in particular, $\2^n$) embeds into a large enough finite partition lattice (see Theorem~413 in \cite{Gr}), and partition lattices are simple (see Theorem~404 in \cite{Gr}), hence by Remark~\ref{remark (2)->(1) simple lattices}, (Ess) holds only for $n\leq 1$.

Now we show that for arbitrary lattices neither (Sub) $\Rightarrow$ (Quo) nor (Quo) $\Rightarrow$ (Sub) holds in general. Using partition lattices again, we can give counterexamples to (Sub) $\Rightarrow$ (Quo).
Since $\2^n$ embeds into a large enough finite partition lattice and partition lattices are simple, we have that $\2^n$ is not a homomorphic image of a partition lattice.

To disprove (Quo) $\Rightarrow$ (Sub), let $ K_n=\mathcal{P}(\{1, \dots, n\}) \cong \2^n$ for some $n\geq4$, and
define a partial order on the set $L:=K_n\times\{0,1\}$ as follows. For
$(a,i),(b,j)\in K_n$, let $(a,i)\leq(b,j)$ iff either $a\leq b$, or $a=b$ and
$i\leq j$. Note that this is the lexicographic order on $L$, and this makes
$L$ a lattice with the following lattice operations (here $\parallel$ stands
for incomparability in $K_n$):%
\[
(a,i)\vee(b,j)=\left\{  \!\!%
\begin{array}
[c]{rl}%
(a,i\vee j), & \text{if }a=b,\\
(a,i), & \text{if }a>b,\\
(b,j), & \text{if }a<b,\\
(a\vee b,0), & \text{if }a\parallel b;
\end{array}
\right.  \qquad\qquad(a,i)\wedge(b,j)=\left\{  \!\!%
\begin{array}
[c]{rl}%
(a,i\wedge j), & \text{if }a=b,\\
(b,j), & \text{if }a>b,\\
(a,i), & \text{if }a<b,\\
(a\wedge b,1), & \text{if }a\parallel b.
\end{array}
\right.
\]
Now $K_n$ is a homomorphic image of $L$ under the homomorphism
$L\rightarrow K_n,~(a,i)\mapsto a$. To see that $K_n$
does not occur as a sublattice of $L$, note that $\{1,2\}$ is a doubly
reducible element in $K_n$, i.e., it can be written as a join as
well as a meet of two incomparable elements: $\{1,2\}=\{1\}\vee
\{2\}=\{1,2,3\}\wedge\{1,2,4\}$. However, there is no doubly reducible element
in $L$, since a nontrivial join in $L$ is always of the form $(a,0)$, and a
nontrivial meet is always of the form $(a,1)$. (It is not necessary to double each element of $K_n$: with a more careful argument, one can construct a counterexample of only 
$2^n+1$ elements.)

Note that the assumption $n\geq4$ was essential in the construction of this
counterexample, since for $n\leq3$, there are no doubly reducible elements in
$K_n$. In fact, one can prove that (Quo) $\Rightarrow$ (Sub) holds for all
lattices (distributive or not) for $n\leq3$ (see Lemma~73 in \cite{Gr}.

Summarizing the results up to this point we know that (Ess) $\Rightarrow$ (Sub) holds, but none of the implications (Sub) $\Rightarrow$ (Ess), (Sub) $\Rightarrow$ (Quo) or (Quo) $\Rightarrow$ (Sub) hold for arbitrary finite lattices. This immediately implies that (Quo) $\Rightarrow$ (Ess) can not hold in general.
The lattice $L=M_3^n$ shows that (Ess) $\Rightarrow$ (Quo) does not hold, either. It is straightforward to verify that the operation $f \in \OL^{(n)}$ defined by 
$$f\bigl((x_{11},\dots,x_{1n}),\dots,(x_{n1},\dots,x_{nn})\bigr)=(x_{11},\dots,x_{n1})$$ commutes with $\lor$ and $\land$ and depends on all of its variables, hence (Ess) holds for $L$.
By Lemma~\ref{lemma CD variety theta=theta_1 x ... x theta_n}, every quotient of $ M_3^n$ is isomorphic to a product of quotients of $M_3$. Since $M_3$ is simple, $L$ only has the quotients $M_3, M_3^2, \dots, M_3^n$, and therefore $\2^n$ does not appear as the quotient algebra of $L$.

 It is also an interesting question to investigate whether any two of the three statements (Ess), (Sub) and (Quo) (of Theorem~\ref{thm essentially n-ary in C* <-> 2^n sublattice <-> 2^n congruence}) imply the third statement in general. First, it is easy to see that (Ess) and (Sub) together do not imply (Quo), but (Ess) and (Quo) imply (Sub), since we have (Ess) $\Rightarrow$ (Sub) and (Ess) $\nRightarrow$ (Quo). We will show that (Sub) and (Quo) together imply (Ess) for arbitrary finite lattices. Let us suppose that (Sub) and (Quo) hold for a finite lattice $L$. Then by (Quo), there is a quotient of $L$ that is isomorphic to $\2^n$, and since $\2$ is a homomorphic image of $\2^n$, we have that $\2$ is a homomorphic image of $L$. Let $\Phi$ denote a surjective homomorphism $\Phi \colon L \rightarrow \2$. Then obviously $\ker(\Phi) \neq L^2$ and by (Sub) we have that $L/\ker(\Phi) \times \dots \times L/\ker(\Phi) = \big( L/\ker(\Phi) \big)^n \cong \2^n$ embeds into $L$. Therefore, by Lemma~\ref{lemma CD variety essentially n-ary in C* <-> A_1/theta_1 x ... x A_n/theta_n <= A}, there exists an essentially $n$-ary operation in $\OL$.

This section gave us some insight to the appearance of $n$-ary operations in the centralizer of an arbitrary finite lattice. We summarize these results in the following proposition.
\begin{proposition}
 Let $L=(L; \lor, \land)$ be an arbitrary finite lattice. Then the following are true:
	\begin{itemize} 
 		\item For any $n$-ary operation $f \in [\lor, \land]^\ast$, there exist unary operations $u_1, \dots, u_n \in [\lor, \land]^\ast$ such that $f(x_1,  \dots, x_n) = u_1(x_1) \lor \dots \lor u_n(x_n)$ and for all $i, j \in \{1, \dots, n\}, i\neq j$ we have $u_i(1) \land u_j(1) = u_1(0)= \dots = u_n(0)$.
		\item If there is an essentially $n$-ary operation in $[\lor, \land]^\ast$, then there is a sublattice of $L$ that is isomorphic to $\2^n$.
		\item If there is a sublattice $L_s$ and a quotient $L_q$ of $L$ that are both isomorphic to $\2^n$, then there is an essentially $n$-ary operation in $[\lor, \land]^\ast$.
	\end{itemize}
\end{proposition}

\section{Centralizer clones over the two-element set} \label{sect:boole}

As promised in the introduction, we are going to determine the centralizer of each clone on $\{0,1\}$ in a fairly simple way.
We use the Post lattice and the notation for Boolean clones from the appendix (see Figure~\ref{fig:post} and Table~\ref{table:boolean clones}). 
First let us record two entirely obvious facts, just for reference:
\begin{fact} \label{fact (C_1 V C_2)* = C_1 * cap C_2 *}
For any clones $C_1, C_2 \leq \OA$ we have $(C_1 \lor C_2)^\ast = C_1^\ast \land C_2^\ast$. (Here $\lor$ and $\land$ denote the join and meet operations of the clone lattice over $A$, i.e., $C_1 \lor C_2 = [ C_1 \cup C_2]$ and $C_1 \land C_2 = C_1 \cap C_2$.) This implies that if $C_1 \leq C_2$ then $C_2^\ast \leq C_1^\ast$.
\end{fact}

\begin{fact} \label{fact 0 in C* <-> C <= Omega_1 ...}
By the definition of the clones $\Omega_0, \Omega_1$ and $S$, for any clone $C \leq \O_{\{0,1\}}$ we have
	\begin{itemize}
		\item $0 \in C^\ast \iff C \leq \Omega_0$,
		\item $1 \in C^\ast \iff C \leq \Omega_1$,
		\item $\neg \in C^\ast \iff C \leq S$.
	\end{itemize}
\end{fact}

We will see that all centralizers over $\{0,1\}$ can be computed using three tools: Theorem~\ref{thm Larose}, Corollary~\ref{cor:log2}, and Proposition~\ref{prop:abelian} below.
This proposition can be found in \cite[Proposition~2.1]{Szend}, but, for the reader's convenience, we include a proof, which is very similar to the proof of Theorem~\ref{thm Larose}.

\begin{proposition} \label{prop:abelian}
Let $A=(A; +)$ be an Abelian group, and let $m(x,y,z)=x-y+z$.
An $n$-ary operation $f \in \OA$  belongs to the centralizer $[m]^\ast$ if and only if there exist unary operations $u_1, \dots, u_n \in [m]^\ast$ such that
\[
f(x_1, \dots, x_n) = u_1(x_1) + \dots + u_n(x_n).
\]
\end{proposition}

\begin{proof}
It is easy to see that the addition commutes with $m$, hence if $f$ is a sum of endomorphisms of $(A;m)$, then $f\in[m]^\ast$. 

Conversely, assume that $f$ is an $n$-ary operation in $[m]^\ast$, and define $u_1, \dots, u_n$ the same way as in the proof of Theorem~\ref{thm Larose}: $u_1(x)=f(x, 0, \dots, 0), \dots, u_n(x)=f(0, \dots, 0, x)$. 
Then we have $u_1, \dots, u_n \in [m]^\ast$, as $f$ and the constant $0$ operation commute with $m$. 
Let us consider the $(n+1)$-ary operation $g(x_1,\dots,x_n,y)=x_1+\dots+x_n-(n-1)y$. The following expression shows that $g\in[m]$:
\begin{align*}
g(x_1,\dots,x_n,y)&=x_1 - y + x_2 - y + x_3 - \dots -y + x_n \\
&= m(\cdots m(m(x_1,y,x_2),y,x_3),\dots,y,x_n).
\end{align*}
(Actually, it is well known and also easy to verify that the elements of $[m]$ are the operations of the form $a_1x_1+\dots+a_nx_n\, (n\in\mathbb{N}^+,a_i\in\mathbb{Z},\sum a_i=1)$, but for the purposes of this proof we only need the operation $g$ above.)
Since $g\in[m]$ and $f\in[m]^\ast$, the operations $f$ and $g$ commute. Applying the definition of commutation to the $(n+1) \times n$ matrix
$$\begin{pmatrix}
x_1 & 0 & 0 & \dots & 0 \\
0 & x_2 & 0 & \dots & 0 \\
\vdots & \vdots & \vdots & \ddots & \vdots \\
0 & 0 & 0 & \dots & x_n \\
0 & 0 & 0 & \dots & 0 
\end{pmatrix},$$
\noindent
we can conclude that
\begin{align*} f(x_1, \dots, x_n)
&= f\big(g(x_1 ,0 ,\dots ,0,0),\ \dots,\ g(0 ,\dots,0 ,x_n ,0)\big) \\
&= g\big(f(x_1, 0, \dots, 0) ,\dots ,f(0, \dots, 0, x_n),f(0,\dots,0)\big)\\
&= u_1(x_1) + \dots + u_n(x_n) - (n-1)\cdot f(0,\dots,0).
\end{align*}
This is almost the required form of $f$; we only need to deal with the constant term  $(n-1)\cdot f(0,\dots,0)$.
However, since $m$ is idempotent, every constant commutes with $m$, thus $u'_n(x_n):=u_n(x_n) - (n-1)\cdot f(0,\dots,0)$ belongs to $[m]^\ast$.
Then we can write $f$ as $f(x_1, \dots, x_n) = u_1(x_1) + \dots + u'_n(x_n)$, and this completes the proof.
\end{proof}

\begin{theorem} \label{thm centralizers on 01}
The centralizers of the clones of Boolean functions are as indicated in Table~\ref{table centralizers on 01} in the appendix.
The clones are grouped by their centralizer clones; the first column shows the 25 primitive positive clones over the two-element set, and the second column lists all clones having the given primitive positive clone as their centralizer.
\end{theorem}
\begin{proof}
Let us recall that a variety $\mathcal{V}$ is congruence distributive if and only if it has, for some $n$, a sequence of terms $J_0(x, y, z), \dots , J_n(x, y, z)$ satisfying the following identities: 
\begin{align*} 
J_0(x, y, z) & = x, \\
J_n(x, y, z) & = z, \\
J_i(x, y, x) & = x \text{ for each } 0 \leq i \leq n, \\
J_i(x, x, y) & = J_{i+1}(x, x, y) \text{ if } i \text{ is even}, \\
J_i(x, y, y) & = J_{i+1}(x, y, y) \text{ if } i \text{ is odd.}
\end{align*} 
These terms $J_i$ are called \emph{Jónnson terms}. If $C \leq \O_{\{0,1\}}$ is a clone, then the existence of a sequence of Jónnson terms in the clone $C$ guarantees that the variety generated by the algebra $(\{0,1\}; C)$ is congruence distributive. By Corollary~\ref{cor:log2}, this implies that $C^\ast \leq \Omega^{(1)}$. 

The clone $SM$ of self-dual monotone Boolean functions is generated by the majority operation $\mu(x,y,z)=xy+xz+yz$, which immediately gives us a sequence of Jónnson terms with $J_0(x,y,z)=x, J_1(x,y,z)=\mu(x,y,z)$ and $J_2(x,y,z)=z$.
We provide a sequence of Jónnson terms in $U^\infty_{01}M$ in Table~\ref{table jonnson terms}; the duals of these operations are Jónsson terms in $W^\infty_{01}M$.
Thus if $C$ contains at least one of the three clones $U^\infty_{01}M, W^\infty_{01}M$ and $SM$ as a subclone, then $C^\ast$ contains only essentially at most unary functions by Corollary~\ref{cor:log2}, and then $C^\ast$ is easy to find using Fact~\ref{fact 0 in C* <-> C <= Omega_1 ...}.
This covers the first six rows of Table~\ref{table centralizers on 01}.

\begin{table}[t] 
\[
\renewcommand{\arraystretch}{1.3}
\begin{tabular}{c c c| c c c c c}
$x$ & $y$ & $z$ & $x=J_0$ & $J_1$ & $J_2$ & $J_3$ & $J_4=z$ \\ \hline
0 & 0 & 0 & 0 & 0 & 0 & 0 & 0 \\
0 & 0 & 1 & 0 & 0 & 0 & 0 & 1 \\
0 & 1 & 0 & 0 & 0 & 0 & 0 & 0 \\
0 & 1 & 1 & 0 & 0 & 0 & 1 & 1 \\
1 & 0 & 0 & 1 & 0 & 0 & 0 & 0 \\
1 & 0 & 1 & 1 & 1 & 1 & 1 & 1 \\
1 & 1 & 0 & 1 & 1 & 0 & 0 & 0 \\
1 & 1 & 1 & 1 & 1 & 1 & 1 & 1 \\
\end{tabular}
\]
\caption{A sequence of Jónnson terms in the clone $U^\infty_{01}M$.} \label{table jonnson terms}
\end{table}

After having determined clones with essentially unary centralizers, there are finitely many clones left to investigate. 
It is easy to see that these clones appear as joins of some of the clones $[0], [1], [\neg],  V_{01}, \Lambda_{01}$ and $L_{01}$. 
According to Fact~\ref{fact (C_1 V C_2)* = C_1 * cap C_2 *}, it suffices to determine the centralizers of these six clones. 
It follows immediately from the definition of $\Omega_0$, $\Omega_1$ and $S$ that $[0]^\ast=\Omega_0$, $[1]^\ast=\Omega_1$ and $[\neg]^\ast=S$.

Theorem~\ref{thm Larose} gives us the centralizer of $V_{01}$: every operation in $V_{01}^\ast$ is of the form $u_1(x_1) \lor \dots \lor u_n(x_n)$, where $u_i(x_i)=x_i$ or $u_i$ is constant for all $i=1,\dots,n$. 
Thus, we have $V_{01}^\ast=[\lor, 0, 1]=V$, and dually, $\Lambda_{01}^\ast=[\land, 0, 1]=\Lambda$. 
Finally, Proposition~\ref{prop:abelian} shows that the centralizer of $L_{01}=[x-y+z]=[x+y+z]$ consists of sums of unary functions, hence $L_{01}^\ast=\{x_1+x_2+\cdots+x_n+c \mid c \in \{0,1\}, n \in \N_0\}=L$.
\end{proof}

The following remark makes it easier to remember the centralizers of all Boolean clones.

\begin{remark}\label{remark partition}
We can group the clones on $\{0,1\}$ by the ``type" of their centralizers. These groups give a partition of the Post lattice into five blocks:
	\begin{itemize}
		\item $\mathcal{C}_{\rm{cd}} := \{ C \leq \O_{\{0,1\}} \mid U^\infty_{01} M \leq C \text { or } W^\infty_{01} M \leq C \text{ or } SM \leq C\}$;
		\item $\mathcal{C}_\lor:= \{V, V_0, V_1, V_{01}\}$;
		\item $\mathcal{C}_\land := \{\Lambda, \Lambda_0, \Lambda_1, \Lambda_{01}\}$;
		\item $\mathcal{C}_{\rm{lin}} := \{L, L_0, L_1, L_{01},SL\}$;
		\item $\mathcal{C}_{\rm{un}} :=\{ C \leq \O_{\{0,1\}} \mid C \leq \Omega^{(1)}\}$.
	\end{itemize}
The ``centralizing" operation $C\mapsto C^\ast$ preserves this partition: for every $C \leq \O_{\{0,1\}}$, we have
	\begin{itemize}
		\item if $C \in \mathcal{C}_{\rm{cd}}$ then $C^\ast \in \mathcal{C}_{\rm{un}}$ (i.e., $C^\ast \leq \Omega^{(1)}$);
		\item if $C \in \mathcal{C}_\lor$ then $C^\ast \in \mathcal{C}_\lor$;
		\item if $C \in \mathcal{C}_\land$ then $C^\ast \in \mathcal{C}_\land$;
		\item if $C \in \mathcal{C}_{\rm{lin}}$ then $C^\ast \in \mathcal{C}_{\rm{lin}}$;
		\item if $C \in \mathcal{C}_{\rm{un}}$, then $C^\ast\geq S_{01}$ (and thus $C \in \mathcal{C}_{\rm{cd}})$.
	\end{itemize}
Figure~\ref{fig:post} in the appendix shows the above partition of the Post lattice  (the five blocks are indicated by different symbols) with primitive positive clones marked by a symbol having an outline. 
Observe that primitive positive clones belonging to the same block have different unary parts most of the time, the only exception being $\Omega_{01}\cap\Omega^{(1)}=S_{01}\cap\Omega^{(1)}=[x]$. 
Thus the observations above together with Fact \ref{fact 0 in C* <-> C <= Omega_1 ...} allow us to find the centralizer of any clone with ease.
\end{remark}

\section*{Acknowledgements}\label{sec:ack}

The authors are grateful to G\'{a}bor Cz\'{e}dli and L\'aszl\'o Z\'adori for helpful discussions.

This research was partially supported by the National Research, Development and Innovation Office of Hungary under grants  KH126581 and K128042, and by the \'{U}NKP-19-3-SZTE-305 new national excellence program and grant TUDFO/47138-1/2019-ITM of the Ministry for Innovation and Technology, Hungary.

\appendix
\section*{Appendix: Clones on the two-element set}

\begin{figure}[h]
\centering
\includegraphics[width=\textwidth]{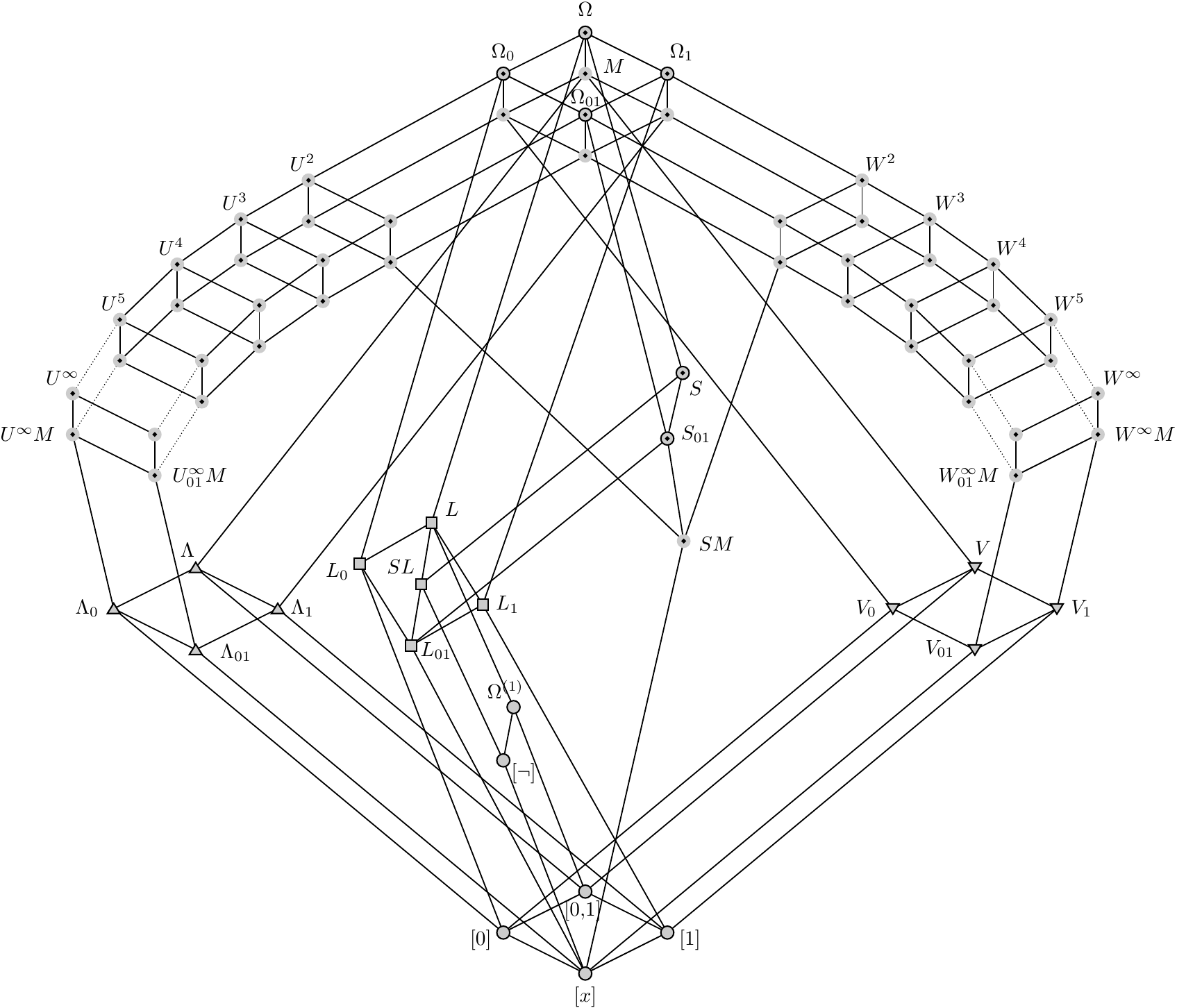}
\caption{The Post lattice.}
\label{fig:post}
\end{figure}

The lattice of all clones on the set $\{0,1\}$ is shown in Figure~\ref{fig:post}.
Different symbols are used according to the partition defined in Remark~\ref{remark partition}; primitive positive clones are indicated by a symbol having an outline, while the gray circles without an outline indicate clones that are not primitive positive.
In Table~\ref{table:boolean clones} we give the definitions of the clones that are labelled on the diagram; the remaining clones can be obtained as intersections of some of these clones.

\vfill

\vfill

\begin{table}[p]
\renewcommand{\arraystretch}{1.6}
\begin{tabular}{l}
\hline
$\Omega$ is the clone of all Boolean functions: $\Omega=\O_{01}$.\\ \hline 
$\Omega_0=\{f \in \Omega \mid f(0,\dots,0)=0\}$ is the clone of 0-preserving functions.\\
$\Omega_1=\{f \in \Omega \mid f(1,\dots,1)=1\}$ is the clone of 1-preserving functions.\\
$\Omega_{01}=\Omega_0 \cap  \Omega_1$ is the clone of idempotent functions.\\ 
In general, if $C$ is a clone, then $C_0=C \cap \Omega_0$, $C_1=C \cap \Omega_1$, and $C_{01}=C_0 \cap C_1$.\\ \hline 
$\Omega^{(1)}$ is the clone of all essentially at most unary functions: $\Omega^{(1)}= [x, \neg x, 0, 1]$.\\ 
$[x]$ is the trivial clone containing only projections.\\ 
$[0]$, $[1]$ and $[0,1]$ are the clones generated by constant operations.\\ 
$[\neg]$ is the clone generated by negation.\\ \hline 
$M=\{ f \in \Omega \mid \x \leq \y \Rightarrow f(\x) \leq f(\y) \}$ is the clone of monotone functions.\\ \hline 
$U^\infty = \{ f \in \Omega^{(n)} \mid n \in \N_0, \exists k \in \{1,\ldots,n\} \colon f(\x)=1 \implies x_k=1 \}$, and\\
$U^\infty M=U^\infty \cap M$, and $U^\infty _{01} M =U^\infty \cap\Omega_{01} \cap M$.\\ \hline 
$W^\infty = \{ f \in \Omega^{(n)} \mid n \in \N_0, \exists k \in \{1,\ldots,n\} \colon f(\x)=0 \implies x_k=0 \}$, and\\ 
$W^\infty M=W^\infty \cap M$ and $W^\infty _{01} M =W^\infty \cap \Omega_{01} \cap M$.\\ \hline 
$S=\{f \in \Omega \mid \neg f(\neg \x)=f(\x)\}$ is the clone of self-dual functions.\\ 
$SM = S \cap M = [\mu]$ where $\mu(x,y,z)$ is the majority function on $\{0,1\}$.\\ \hline 
$\Lambda=\{ x_1 \land \cdots \land x_n \mid n \in \N^+\} \cup [0,1]=[\land,0,1]$.\\ 
$\Lambda_0=\Lambda \cap \Omega_0=\{ x_1 \land \cdots \land x_n \mid n \in \N^+\} \cup [0]=[\land,0]$.\\ 
$\Lambda_1=\Lambda \cap \Omega_1=\{ x_1 \land \cdots \land x_n \mid n \in \N^+\} \cup [1]=[\land,1]$.\\ 
$\Lambda_{01}=\Lambda \cap \Omega_{01}=\{ x_1 \land \cdots \land x_n \mid n \in \N^+\}=[\land]$.\\ \hline 
$V=\{ x_1 \lor \cdots \lor x_n \mid n \in \N^+\} \cup [0,1]=[\lor,0,1]$.\\ 
$V_0=V \cap \Omega_0=\{ x_1 \lor \cdots \lor x_n \mid n \in \N^+\} \cup [0]=[\lor,0]$.\\ 
$V_1=V \cap \Omega_1=\{ x_1 \lor \cdots \lor x_n \mid n \in \N^+\} \cup [1]=[\lor,1]$.\\ 
$V_{01}=V \cap \Omega_{01}=\{ x_1 \lor \cdots \lor x_n \mid n \in \N^+\}=[\lor]$.\\ \hline 
$L=\{x_1+\cdots+x_n+c \mid c \in \{0,1\}, n \in \N_0\}=[x+y,1]$.\\ 
$L_0=L \cap \Omega_0=\{x_1+\cdots+x_n \mid n \in \N_0\}=[x+y]$.\\ 
$L_1=L \cap \Omega_1=\{x_1+\cdots+x_n+(n+1\bmod 2) \mid n \in \N_0\}=[x+y+z,1]$.\\ 
$L_{01}=L \cap \Omega_{01}=\{x_1+\cdots+x_n \mid n \text{ is odd}\}=[x+y+z]$.\\ 
$SL=S \cap L=\big\{x_1+\cdots+x_n+c \mid n\text{ is odd and }c\in\{0,1\}\big\}=[x+y+z,x+1]$.\\ \hline
\end{tabular}
\bigskip
\caption{Definitions of some clones of Boolean functions.}
\label{table:boolean clones}
\end{table}

\floatstyle{ruled}
\begin{table}[p]
\renewcommand{\arraystretch}{1.9}
\begin{tabular}{|p{0,6cm}|p{11,6cm}|}  \hline
$P$ 		&  all clones $C \leq \O_{\{0,1\}}$ such that $C^\ast = P$ \\ \Xhline{3\arrayrulewidth}
$[x]$ 		& $\Omega, M$\\ \hline
$[0]$		& $\Omega_0, M_0, U^k, U^\infty, U^{k}M, U^{\infty}M$ (for all $k \in \N^+$)\\ \hline
$[1]$ 		& ${\Omega_1}, {M_1}, W^k, W^\infty, W^k M, W^\infty M$ (for all $k \in \N^+ $)\\ \hline
$[0,1]$ 	& $\Omega_{01}, M_{01}, U_{01}^k, U_{01}^\infty, U_{01}^k M, U_{01}^\infty M, W_{01}^k, W_{01}^\infty, W_{01}^k M, W_{01}^\infty M $ (for all $k \in \N^+$)\\ \hline
$[\neg]$ 			& $S$\\ \hline
$\Omega^{(1)}$ 		& $S_{01}, SM$\\ \hline
$L_{01}$ 	& $L$\\ \hline
$L_0$ 	& $L_0$\\ \hline
$L_1$ 	& $L_1$\\ \hline
$L$ 		& $L_{01}$\\ \hline
$SL$ 		& $SL$\\ \hline
$\Lambda_{01}$ 		& $\Lambda$\\ \hline
$\Lambda_0$ 		& $\Lambda_0$\\ \hline
$\Lambda_1$ 		& $\Lambda_1$\\ \hline
$\Lambda$ 			& $\Lambda_{01}$\\ \hline
$V_{01}$	& $V$\\ \hline
$V_0$		& $V_0$\\ \hline
$V_1$ 	& $V_1$\\ \hline
$V$		& $V_{01}$\\ \hline
$S_{01}$ 	& $\Omega^{(1)}$\\ \hline
$S$		& $[\neg]$\\ \hline
$\Omega_{01}$ 		& $[0,1]$\\ \hline
$\Omega_0$			& $[0]$\\ \hline
$\Omega_1$			& $[1]$\\ \hline
$\Omega$ 			& $[x]$\\ \hline
\end{tabular}
\bigskip
\caption{The centralizers of all clones of Boolean functions.}\label{table centralizers on 01}
\end{table}

\end{document}